\documentclass[a4paper,12pt]{amsart}
\usepackage{amsfonts}
\usepackage{amssymb}
\usepackage{ifthen}
\usepackage{graphicx}
\nonstopmode \numberwithin{equation}{section}
\usepackage[usenames]{color}
\newtheorem{thm}{Theorem}[section]
\newtheorem{lem}{Lemma}[section]
\newtheorem{cor}{Corollary}[section]

\newtheorem{cl}{Claim}[section]
\newtheorem{ca}{Case}[section]
\newtheorem{sca}{Subcase}[section]
\newtheorem{scl}[section]{Subclaim}
\newtheorem{conj}[equation]{Conjecture}

\theoremstyle{definition}
\newtheorem{defn}{Definition}[section]
\newtheorem{op}[equation]{Open Problem}
\newtheorem{ques}[equation]{Question}
\newtheorem{rem}{Remark}[section]
\newtheorem{exam}[equation]{Example}

\newcounter {own}
\def\theown {\thesection       .\arabic{own}}

\newenvironment{pf}[1][]{%
 \vskip 3mm
 \noindent
 \ifthenelse{\equal{#1}{}}%
  {{\slshape Proof. }}%
  {{\slshape #1.} }%
 }%
{\qed\bigskip}

\newcounter{alphabet}
\newcounter{tmp}
\newenvironment{Thm}[1][]{\refstepcounter{alphabet}%
\bigskip%
\noindent%
{\bf Theorem \Alph{alphabet}}%
\ifthenelse{\equal{#1}{}}{}{ (#1)}%
{\bf .} \itshape}{\vskip 8pt}

\makeatletter
\newcommand{\Ref}[1]{\@ifundefined{r@#1}{}{\setcounter{tmp}{\ref{#1}}\Alph{tmp}}}
\makeatother

\newenvironment{Lem}[1][]{\refstepcounter{alphabet}%
\bigskip%
\noindent%
{\bf Lemma \Alph{alphabet}}%
{\bf .} \itshape}{\vskip 8pt}

\newcommand{\diam}{{\operatorname{diam}}}

\newcommand{\dist}{{\operatorname{dist}}}

\def\be{\begin{equation}}
\def\ee{\end{equation}}

\newcommand{\ben}{\begin{enumerate}}
\newcommand{\een}{\end{enumerate}}

\newcommand{\blem}{\begin{lem}}
\newcommand{\elem}{\end{lem}}
\newcommand{\bthm}{\begin{thm}}
\newcommand{\ethm}{\end{thm}}
\newcommand{\bcor}{\begin{cor}}
\newcommand{\ecor}{\end{cor}}
\newcommand{\beg}{\begin{exam}}
\newcommand{\eeg}{\end{exam}}
\newcommand{\begs}{\begin{examples}}
\newcommand{\eegs}{\end{examples}}
\newcommand{\bdefe}{\begin{defn}}
\newcommand{\edefe}{\end{defn}}
\newcommand{\bprob}{\begin{prob}}
\newcommand{\eprob}{\end{prob}}
\newcommand{\bques}{\begin{ques}}
\newcommand{\eques}{\end{ques}}
\newcommand{\bei}{\begin{itemize}}
\newcommand{\eei}{\end{itemize}}
\newcommand{\bcon}{\begin{conj}}
\newcommand{\econ}{\end{conj}}
\newcommand{\bop}{\begin{op}}
\newcommand{\eop}{\end{op}}

\newcommand{\bca}{\begin{ca}}
\newcommand{\eca}{\end{ca}}
\newcommand{\bsca}{\begin{sca}}
\newcommand{\esca}{\end{sca}}

\newcommand{\bcl}{\begin{cl}}
\newcommand{\ecl}{\end{cl}}

\newcommand{\bscl}{\begin{scl}}
\newcommand{\escl}{\end{scl}}

\newcommand{\bcons}{\begin{conjs}}
\newcommand{\econs}{\end{conjs}}
\newcommand{\bprop}{\begin{propo}}
\newcommand{\eprop}{\end{propo}}
\newcommand{\br}{\begin{rem}}
\newcommand{\er}{\end{rem}}
\newcommand{\brs}{\begin{rems}}
\newcommand{\ers}{\end{rems}}
\newcommand{\bo}{\begin{obser}}
\newcommand{\eo}{\end{obser}}
\newcommand{\bos}{\begin{obsers}}
\newcommand{\eos}{\end{obsers}}
\newcommand{\bpf}{\begin{pf}}
\newcommand{\epf}{\end{pf}}
\newcommand{\ba}{\begin{array}}
\newcommand{\ea}{\end{array}}
\newcommand{\beq}{\begin{eqnarray}}
\newcommand{\beqq}{\begin{eqnarray*}}
\newcommand{\eeq}{\end{eqnarray}}
\newcommand{\eeqq}{\end{eqnarray*}}

\newcounter{minutes}
\divide\time by 60
\newcounter{hours}
\multiply\time by 60 \addtocounter{minutes}{-\time}

\begin{document}

\bibliographystyle{amsplain}
\title{Sphericalization and flattening with their applications in quasimetric measure spaces}

\author{Qingshan Zhou}
\address{Qingshan Zhou, School of Mathematics and Big Data, Foshan University,  Foshan, Guangdong 528000, People's Republic
of China} \email{q476308142@qq.com}

\author{Yaxiang Li}
\address{Yaxiang Li,  Department of Mathematics, Hunan First Normal University, Changsha,
Hunan 410205, P.R.China}
\email{yaxiangli@163.com}

\author{Xining Li${}^{~\mathbf{*}}$}
\address{Xining Li,  Sun Yat-sen University, Department of Mathematics, Guangzhou 510275, People's Republic
of China} \email{lixining3@mail.sysu.edu.cn}

\def\thefootnote{}
\footnotetext{\date{\today}
} \makeatletter\def\thefootnote{\@arabic\c@footnote}\makeatother

 \keywords{
Sphericalization, flattening,
doubling condition, Ahlfors regular, Loewner condition, quasi-metric spaces, quasi-m\"{o}ius, Gromov hyperbolic spaces.\\
${}^{\mathbf{*}}$ Corresponding author}

\begin{abstract} The main purpose of the note is to explore the invariant properties of sphericalization and flattening and their applications in quasi-metric spaces. We show that sphericalization and flattening procedures on a quasimetric spaces preserving properties such as  Ahlfors regular  and doubling property. By using these properties, we generalize a recent result in \cite{WZ}. We also show that the Loewner condition can be  preserved under quasim\"obius   mapping between two $Q$-Ahlfors regular spaces. Finally, we prove that the $Q$-regularity of $Q$-dimensional Hausdorff measure  of Bourdon metric are coincided with Hausdorff measure of    Hamenst\"adt metric defined on the boundary at infinity of a Gromov hyperbolic space.

\end{abstract}

\thanks{The first author was supported by NNSF of
	China (No. 11571216), the second author   was   supported by NNSF of
	China (Nos. 11601529,  11671127)  and the third author was supported by  NNSF of China (No. 11701582).}

\maketitle{} \pagestyle{myheadings} \markboth{}{}

\section{Introduction and main results}\label{sec-1}

The original idea of {\it sphericalization} and {\it flattening} comes from the work of Bonk and Kleiner \cite{BK02} in defining a metric on the one point compactification of an unbounded locally compact metric space. The first class of deformation, sphericalization, is a generalization of the deformation from the Euclidean distance on $\mathbb{R}^n$ to the chordal distance on $\mathbb{S}^n$. The second class of flattening deformation is a generalization of inversion on punctured $\mathbb{S}^n$.

It was shown in \cite{BHX} that these two conformal transformations are dual in the sense that if one starts from a bounded metric space, then performs a flattening transformation followed by a sphericalization, then the object space is bilipschitz equivalent to the original space. This duality comes from the idea that the stereographic projection between the Euclidean space and the Riemann sphere can be realized as a special case of inversion. Sphericalization and flattening have a lot of applications in the area of geometric function theory and asymptotic geometry, such as \cite{BB,BHX,BuSc,HSX}.

By using the sphericalization (named by {a warping process} in \cite{W}), Wildrick  obtained the quasisymmetric parameter of an unbounded $2$-dimensional metric planes. In \cite{JJ}, Jordi proved that two visual geodesic Gromov hyperbolic spaces are roughly quasi-isometric if and only if their boundary at infinity are power quasim\"obius equivalent by virtue of the flattening deformation. Moreover, Mineyev \cite{M} studied the metric conformal structures on the idea boundaries of hyperbolic complexes via sphericalization.  Durand-Cartagena, the third author, Shanmugalingam \cite{DL1,DL2,LS} explored the preservation of Ahlfors, doubling measure and Poincar\'e inequality under sphericalization and flattening within different assumptions.

The main purpose of this paper is to investigate the invariant properties of sphericalization and flattening and their applications in quasi-metric measure spaces. There are many materials for quasi-metric spaces and related geometric function theory in this setting, such as a quasi-metric space quasisymmetric embedding onto Euclidean spaces and $H^p$ theory and Lipschitz function theory on  spaces of homogeneous types. We exhibit some examples of quasimetric spaces which are mainly motivated  from this setting.

$(1)$ $Z=\mathbb{R}^n,\rho(x,y)=\sum_{i=1}^n |x_i-y_i|^{\alpha_i}$, where $\alpha_1,\alpha_2,...,\alpha_n$ are positive numbers, not all equal. In general, one can easily see that  the  nonnegative symmetric function is  not necessary a metric but a quasimetric. This example follows from  Coifman and Weiss \cite[\S2(2)]{CW} and they have referred to such  property as {\it nonisotropy}.

$(2)$ An often mentioned example is the {\it snowflake} transformation $d\to d^{\varepsilon}$ of a metric which induces a quasimetric space. The snowflake transformation play an important role in the Assouad embedding theorem \cite{Ass}.

$(3)$ Coifman and Weiss \cite{CW} have given the precise definition for {\it a space of homogeneous type} which is a quasi-metric space  carrying a doubling measure and developed the basic theory of Hardy spaces in this setting. A number of classical examples have also been presented by them such as  compact Riemannian manifolds with natural distances and measure,  boundary of  smooth and bounded pseudo-convex  domains in $\mathbb{C}^n$ with  non-isotropic quasimetric and Lebesgue surface measure, see \cite[\S2(2)]{CW} for more details and several examples.

$(4)$ The idea of {\it deforming metric spaces} by doubling measures is due to David and Semmes, see \cite{DS} and \cite{S}.  Every uniformly perfect homogeneous space can be deformed into Ahlfors regular space, which is quasisymmetric equivalent to the original space through the identity map.

$(5)$ It is not difficult to see that quasi-metrics induced by the Gromov products can be used to define a canonical quasiconformal gauge on the boundary at infinity of hyperbolic spaces in the sense of Gromov, see \cite{BuSc}. Recently many researchers are interested in the interplay between interior and boundary  in quasiconformal geometry for Gromov hyperbolic spaces, see \cite{BHX,BuSc,JJ}.

$(6)$ Buckley, Herron and Xie extended the classical {\it inversion} or reflection about the unit sphere centered at the origin to metric space $(X,d)$ in \cite{BHX}. Let $p\in X$ be given. We find that
$$i_p(x,y)=\frac{d(x,y)}{d(x,p)d(y,p)}$$
defines a quasi-metric on $X\setminus\{p\}$ and the identity map $(X\setminus\{p\},d)\to (X\setminus\{p\},i_p)$ is quasim\"{o}bius. Further, Buyalo and Schroeder introduced a more general $\lambda$-inversion with this idea via an admissible function $\lambda$ on arbitrary quasi-metric spaces in \cite{BuSc}.

In \cite{LS}, the third author and Shanmugalingam showed that the process of sphericalization and flattening preserved the Ahlfors regular and doubling measure in metric spaces. In this paper, we first generalize the work in \cite{LS} to quasimetric spaces. We get the preservation of Ahlfors regularity of a quasimetric space under these two transformations, see Theorems \ref{z-1} and \ref{z-2}. And by using these results and David-Semmes's conformal deformation \cite{DS}, we prove the invariance of doubling measure on uniformly perfect quasimetric spaces under sphericalization and flattening, see Theorem \ref{z-6}. We note that this also provide a new proof for the corresponding results in metric spaces in \cite{LS}.
As a direct application of the above results, we improve the recent work of \cite{WZ} as follows.

\begin{thm}\label{m-thm-1} A weakly quasim\"obius map from a uniformly perfect doubling quasimetric space to a quasimetric space is quasim\"obius if and only if the image space is uniformly perfect and doubling.
\end{thm}

Our next motivation comes from Heinonen and Koskela's celebrated work on the equivalence of quasiconformality and quasisymmetry between two metric spaces in \cite{HeK}. They introduced the concept of {\it Loewner spaces}, which has many applications in studying Sobolev spaces,  quasiconformal theory in metric spaces. It should be noted that Tyson \cite{Tys} answered positively to a conjecture proposed by Heinonen and Koskela \cite[Section 8.7]{HeK} in proving that the $Q$-Loewner condition can be preserved under quasisymmetric maps between two $Q$-regular spaces. Thus it is  natural to ask whether  the Loewner condition is preserved under sphericalization and flattening (more general quasim\"obius mappings). We obtain the following theorem:

\begin{thm}\label{m-thm-2}
Let $(X,d,\mu)$ and $(Y,\sigma,\nu)$ be locally compact $Q$-regular metric measure spaces with $Q>1$. If $f:X\to Y$ is quasi-m\"{o}bius and $X$ is $Q$-Loewner, then $Y$ is also $Q$-Loewner.
\end{thm}
\br

 It is worth to mention that the third author and Shanmugalingam studied the invariance of Poincar\'e inequality under sphericalization and flattening transformations in their recent work \cite{LS}. This is as well one of the main motivation of the current paper. We shall explain the connection between Theorem \ref{m-thm-2} and Li-Shanmugalingam's results in \cite{LS}.

First, Heinonen and Koskela demonstrated that $Q$-Loewner condition and  $(1,Q)$-Poincar\'e inequality are equivalent in proper, Ahlfors $Q$-regular and $\varphi$-convex metric measure spaces, see \cite[ Corollary 5.13]{HeK}. In \cite[Theorem 1.1]{LS}, it was shown by the third author and Shanmugalingam that for a complete doubling metric measure space which  admits a $(1,p)$-Poincar\'e inequality with $1\leq p<\infty$, if in addition the sphericalizad (or flattened) space is annular quasiconvex, then the deformed space also admits a $(1,p)$-Poincar\'e inequality. Notice that we do not need any extra assumptions concerning the connectivity or completeness of the spaces in Theorem \ref{m-thm-2}. Moreover, the Ahlfors $Q$-regularity is not used in showing  the invariance  of $Q$-Loewner condition under  sphericalization and flattening. It is just adapted to the $Q$-Loewner condition because each $Q$-Loewner space satisfies a lower mass  estimate, see \cite[Theorem 3.6]{HeK}. By using a deformed cross-ratio introduced in \cite{BK}, our proof is direct and simple. For more backgrounds in this line see \cite{BHX,HeK,KZ,Kor} and the references therein.

\er

Finally, we give an application which concerns the conformal transformation to the boundary at infinity of Gromov hyperbolic spaces. We investigate the interplay between the Hausdorff measure with respect to the Bourdon metrics and Hamenst\"adt metrisc on the boundary at infinity of a Gromov hyperbolic space, which states as follows.

\begin{thm}\label{m-thm-3}
 Let $X$ be a Gromov hyperbolic metric space and $\partial_\infty X$ the boundary at infinity of $X$, and let $Q>0$. Then $(\partial_\infty X, d_B)$ is Ahlfors $Q$-regular if and only if the space $(\partial_\infty X,d_H)$ is Ahlfors $Q$-regular with respect to their $Q$-dimensional Hausdorff measure, where  $d_B$ and $d_H$ denote the Bourdon metrics and Hamenst\"adt metrics with the same parameter, resepectively.
\end{thm}

  The presentation of this paper is organized as follows. In Section 2, we review some backgrounds and notions of the paper. And then we discuss the preservation of Ahlfors regular quasimetric spaces under sphericalization and flattening in Section 3. After that, we consider the preservation of doubling measure in section 4. Finally, in the last section, we collect the previous results to prove Theorem 1.1, Theorem 1.2 and Theorem 1.3.

%

\section{Preliminary and Notations}

\subsection{Quasi-metric space miscellanea}
To begin our discussion, it is convenient to introduce the concept of a quasi-metric space. Let $K\geq 1$. A {\it quasi-metric} on a set $X$ is a function $\rho:X\times X\to [0,\infty)$ that is symmetric, and $\rho(x,y)=0$ if and only if $x=y$, and  satisfies the following
$$\rho(x,z)\leq K(\rho(x,y)\vee\rho(y,z))$$
for all $x,$ $y,$ $z\in X$. Then a quasi-metric space is a set $X$ together with a quasi-metric. Here and hereafter, we use the notations:
$$r\vee s=\max\{r,s\}\;\;\;\;\;\; \mbox{and}\;\;\;\;\;\; r\wedge s=\min\{r,s\}$$
for $r,$ $s\in \mathbb{R}$.
The quasi-metric balls ($\rho$-balls) are denoted as $B(x,r)=B_\rho(x,r)=\{y\in X: \rho(x,y)<r\}$. Clearly, a metric space is a $K$-quasimetric space with $K=2.$

A quasi-metric space $(X, \rho)$ satisfies a {\it doubling condition} if there is a constant $C$ such that every (quasi-metric) ball $B$ in $X$ can be covered by at most $C$ balls of half the radius of $B$. A positive Borel measure $\mu$ on a quasi-metric space $(X,\rho)$ is a {\it doubling measure} if there is a constant $C_\mu$ such that
$$\mu(2B)\leq C_\mu\mu(B)$$
for all balls $B$. We also record the following result for later use.

\begin{Thm}\label{z-8}$($\cite[Theorem $13.3$]{Hei}$)$
A complete doubling quasi-metric space carries a doubling measure.
\end{Thm}

A quasi-metric space $(X,\rho)$ is called {\it uniformly perfect}, if there is a constant $\tau\in (0,1)$, such that for each $x\in X$ and every $r>0$ for which the set $X\setminus B(x,r)$ is nonempty, we have that $B(x,r)\setminus B(x, \tau r)$ is nonempty.

A quasi-metric space $(X, \rho)$ is said to be {\it Ahlfors $Q$-regular} if $(X, \rho)$ admits a positive Borel measure $\mu$ such that
\be\label{n-1} C^{-1}R^Q\leq \mu(B(x,R)) \leq C R^Q\ee
for all $x\in X$ and $0<R< \diam_\rho(X)$ (It is possible that $\diam_\rho(X)=\infty$), where $C\geq 1$ and $Q>0$ are constants. Note that Ahlfors regular spaces are necessarily doubling. For instance, the Euclidean space $\mathbb{R}^n$ with Lebesgue measure satisfies the Ahlfors $n$-regularity.

\subsection{Modulus of curves and Loewner spaces}
In order to introduce the definition of Loewner spaces, we have to mention about the modulus of curve family. Let $(X,d,\mu)$ be a metric measure space, let $\Gamma$ be a family of nonconstant curves in $X$. We say a Borel function $g$ to be {\it admissible} if
$$\int_{\gamma}gds\ge 1$$
for all locally rectifiable curves $\gamma\in\Gamma$. Let $Q\ge 1$, the {\it $Q$-modulus} of $\Gamma$ is denoted as
$$mod_Q \Gamma=\inf \int_X g^Q d\mu$$
for all admissible function of $\Gamma$. Given a pair of sets $E,F$ in $X$, we denote $\Gamma(E,F)$  as the collection of  rectifiable curve connecting $E$ and $F$ in $X$.

Let $(X,d,\mu)$ be an Ahlfors $Q$-regular metric measure space. Given two disjoint close continua $E,F \subset X$. We denote the {\it relative separation} as follows
$$\Delta(E,F)=\frac{d(E,F)}{ \diam E\wedge\diam F}.$$
We say $(X,d,\mu)$ is {\it $Q$-Loewner} if there exists a homeomorphism $\eta:[0,\infty)\rightarrow [0,\infty)$ such that for each pair of disjoint continua $E,F$ in $X$, it satisfies the following inequality:
$$mod_Q(\Gamma(E,F))\ge \eta(\Delta(E,F)).
$$

\subsection{Mappings on quasi-metric spaces}

For a tetrad $a,$ $b,$ $c,$ $d$ in a quasi-metric space $(Z,\rho)$, its {\it cross ratio} is defined by the number
$$r(a,b,c,d)=\frac{\rho(a,c)\rho(b,d)}{\rho(a,b)\rho(c,d)} \,\,\text{with}\,a \neq b,c\neq d.$$

Suppose that $\eta$ and $\theta$ are homeomorphisms from $[0, \infty)$ to $[0, \infty)$, and that $f:$ $(Z_1,\rho_1)\to (Z_2,\rho_2)$ is an embedding between two quasi-metric spaces. Then  we call $f$ is {\it $L$-bilipschitz} for some $L\geq 1$ if
$$L^{-1}\rho_1(x,y)\le \rho_2(f(x),f(y))\le L\rho_1(x,y)$$
for all $x,y\in Z_1.$ Also, $f$ is said to be {\it $\eta$-quasisymmetric} if
$$\rho_1(x,a)\leq t \rho_1(x,b)\;\;\;\;\;\;\;\; \mbox{implies  that}\;\;\;\;\;\;\;\; \rho_2(f(x),f(a))\leq \eta(t) \rho_2(f(x),f(b))$$
for all $x$, $a$ and $b$ in $Z_1$. On the other hand, $f$ is called
{\it $\theta$-quasim\"obius} if
$$r(a,b,c,d)\leq t\;\;\;\;\;\;\;\; \mbox{implies  that}\;\;\;\;\;\;\;\; r(f(a),f(b),f(c),f(d))\leq \theta(t)$$
for all $a,$ $b,$ $c,$ $d$ in $Z_1$. In particular, if $\theta(t)=t$, then $f$ is a M\"{o}bius mapping. 

 We review the concept of {\it weakly quasim\"obius mapping} here in order to discuss the statement of Theorem \ref{m-thm-1}. We say  $f$ is {\it $(h,H)$-weakly quasim\"obius} if there exist constants $h>0$ and $H\geq 1$ such that
$$r(a,b,c,d)\leq h\;\;\;\;\;\;\;\; \mbox{implies  that}\;\;\;\;\;\; r(f(a),f(b),f(c),f(d))\leq H$$
for all $a,$ $b,$ $c,$ $d$ in $Z_1$.
 In general, it is not easy to determine whether a homeomorphism is quasim\"obius or not. In spirit of this consideration, the authors in  \cite{WZ} introduced this concept and systematically investigated the relationships between weakly quasim\"obius  and quasim\"obius mapping on doubling quasi-metric measure spaces.
\medskip

The following condition was introduced by V\"ais\"al\"a \cite{Vai-5} to investigate the relation between quasim\"{o}bius and quasisymmetric maps. Suppose that both $(Z_1,\rho_1)$ and $(Z_2,\rho_2)$ are bounded quasi-metric spaces. A homeomorphism $f:$ $(Z_1,\rho_1)\to (Z_2,\rho_2)$ is said to satisfy the {\it  $\lambda$-three-point condition} if there are constant $\lambda\geq 1$ and triad $z_1,$ $z_2,$ $z_3$ in $(Z_1, \rho_1)$ such that
$$\rho_1(z_i,z_j)\geq \frac{1}{\lambda}\diam(Z_1)\;\;\mbox{and}\;\; \rho_2(f(z_i),f(z_j))\geq \frac{1}{\lambda}\diam(Z_2)$$
for all $i\neq j\in\{1,2,3\}$. With the help of this condition, he proved

\begin{Thm}\label{z-4}$($\cite{Vai-5}$)$
Suppose that $(Z_1,\rho_1)$ and $(Z_2,\rho_2)$ are bounded quasi-metric spaces and that $f:(Z_1,\rho_1) \to (Z_2,\rho_2)$ satisfies the $\lambda$-three point condition. Then $f$ is quasim\"{o}bius if and only if it is  quasisymmetric, quantitatively.
\end{Thm}

\subsection{Gromov hyperbolic spaces}
Following the terminology of \cite{BuSc}, we say that $(X,d)$ is a {\it Gromov hyperbolic} metric space, if there is a constant $\delta\geq 0$ such that
$$(x|y)_w\geq \min\{(x|z)_w,(z|y)_w\}-\delta$$
for all $x,y,z,w\in X$, where $(x|y)_w$ denotes the {\it Gromov product} with respect to $w$ defined by
$$(x|y)_w=\frac{1}{2}[d(x,w)+d(y,w)-d(x,y)].$$

\bdefe Suppose $(X, d)$ is a Gromov $\delta$-hyperbolic metric space for some constant $\delta\geq 0$.\ben
\item
A sequence $\{x_i\}$ in $X$ is called a {\it Gromov sequence} if $(x_i|x_j)_w\rightarrow \infty$ as $i,$ $j\rightarrow \infty.$
\item
Two Gromov sequences $\{x_i\}$ and $\{y_j\}$ are said to be {\it equivalent} if $(x_i|y_i)_w\rightarrow \infty$.
\item
The {\it Gromov boundary} or {\it the boundary at infinity} $\partial_\infty X$ of $X$ is defined to be the set of all equivalent classes.
\item
For $a\in X$ and $\eta\in \partial_\infty X$, the Gromov product $(a|\eta)_w$ of $a$ and $\eta$ is defined by
$$(a|\eta)_w= \inf \big\{ \liminf_{i\rightarrow \infty}(a|b_i)_w:\; \{b_i\}\in \eta\big\}.$$
\item
For $\xi,$ $\eta\in \partial_\infty X$, the Gromov product $(\xi|\eta)_w$ of $\xi$ and $\eta$ is defined by
$$(\xi|\eta)_w= \inf \big\{ \liminf_{i\rightarrow \infty}(a_i|b_i)_w:\; \{a_i\}\in \xi\;\;{\rm and}\;\; \{b_i\}\in \eta\big\}.$$
\een
\edefe

For $0<\varepsilon<\min\{1,\frac{1}{5\delta}\}$, define
$$\rho_{w,\varepsilon}(\xi,\zeta)=e^{-\varepsilon(\xi|\zeta)_w}$$
for all $\xi,\zeta$ in the Gromov boundary of $X$ with convention $e^{-\infty}=0$.

We now define
$$d_\varepsilon(\xi,\zeta):=\inf\{\sum_{i=1}^{n} \rho_\varepsilon (\xi_{i-1},\xi_i):n\geq 1,\xi=\xi_0,\xi_1,...,\xi_n=\zeta\in \partial_\infty X\}.$$
Then $(X,d_{w,\varepsilon})$ is a metric space with
\be\label{l-0.1}\frac{1}{2}\rho_{w,\varepsilon}\leq d_{w,\varepsilon}\leq \rho_{w,\varepsilon},\ee
and we call $d_{w,\varepsilon}$ to be the {\it Bourdon metric}  \cite[p 22]{BuSc}
 of $\partial_\infty X$ based at $w$ with parameter $\varepsilon$.

Let $\xi\in \partial_\infty X$. We say that $b:X\to \mathbb{R}$ is a {\it Busemann function} based on $\xi$, denoted by $b=b_{\xi,w}\in \mathcal{B}(\xi)$, if for every $w\in X$, we have
$$b(x)=b_{\xi,w}(x)=b_\xi(x,w)=(\xi|w)_x-(\xi|x)_w\;\;\;\;\;\;\;\;\;\;\;\;\mbox{for}\;x\in X.$$

We next define the Gromov product of $x,y\in X$ based at the Busemann function $b=b_{\xi,w}\in \mathcal{B}(\xi)$ by
$$(x|y)_b=\frac{1}{2}(b(x)+b(y)-d(x,y)).$$
Similarly, for $x\in X$ and $\eta\in \partial_\infty X$, the Gromov product $(x|\eta)_b$ of $x$ and $\eta$ is defined by
$$(x|\eta)_b= \inf \big\{ \liminf_{i\rightarrow \infty}(x|z_i)_b:\; \{z_i\}\in \eta\big\}.$$
For points $(\xi_1,\xi_2)\in \partial_\infty X\times  \partial_\infty X  \setminus (\xi,{\xi})$, we define their Gromov product based at $b$ by
$$(\xi_1|\xi_2)_b=\inf\big\{\liminf_{i\to\infty} (x_i|y_i)_b: \{x_i\}\in\xi_1 , \{y_i\}\in\xi_2\}.$$

Next, we review the definition of {\it Hamenst\"adt metric} of $\partial_\infty X$ based at $\xi$ or a Busemann function $b=b_{\xi,w}\in \mathcal{B}(\xi)$. For $\varepsilon>0$ with $e^{22\varepsilon\delta}\leq 2$, define
$$\rho_{b,\varepsilon}(\xi_1,\xi_2)= e^{-\varepsilon(\xi_1|\xi_2)_b}\;\;\;\;\;\;\;\mbox{for all}\;\xi_1,\xi_2\in \partial_\infty X.$$
Then for $i=1,2,3$ with $\xi_i\in \partial_\infty X$, we have
$$\rho_{b,\varepsilon}(\xi_1,\xi_2)\leq e^{22\varepsilon\delta} \max\{\rho_{b,\varepsilon}(\xi_1,\xi_3),\rho_{b,\varepsilon}(\xi_3,\xi_2)\}.$$
That is, $\rho_{b,\varepsilon}$ is a $K'$-quasi-metric on $X$ with $K'=e^{22\varepsilon\delta}\leq 2$. We now define
$$\sigma_{b,\varepsilon}(x,y):=\inf\{\sum_{i=1}^{n} \rho_{b,\varepsilon} (x_{i-1},x_i):n\geq 1,x=x_0,x_1,...,x_n=y\in X\}.$$
Again by \cite[Lemma $2.2.5$]{BuSc}, $(X, \sigma_{b,\varepsilon})$ is a metric space with
\be\label{l-0.2}\frac{1}{2}\rho_{b,\varepsilon} \leq \sigma_{b,\varepsilon}\leq \rho_{b,\varepsilon}.\ee

Then $\sigma_{b,\varepsilon}$ is called a {\it Hamenst\"adt metric} on $\partial_\infty X$ based at $\xi$ or the Busemann function $b$ with parameter $\varepsilon$.

\section{Transformations of quasimetric measure spaces}
\subsection{Sphericalization and flattening of quasi-metric structures}
Given an unbounded quasi-metric space $(X,\rho)$ and a base point $a\in X$, we consider the one-point extension $\dot{X}=X\cup \{\infty\}$ and define the density function $\rho_a:\dot{X}\times \dot{X}\to[0,\infty)$ as follows
\be\label{q-1} \rho_a(x,y)=\rho_a(y,x)=\begin{cases}
\displaystyle\; \frac{\rho(x,y)}{[1+\rho(x,a)][1+\rho(y,a)]},\;
\;\;\;\;\mbox{if}\;\;x,y\in X,\\
\displaystyle\;\;\;\;\;\frac{1}{1+\rho(x,a)},\;\;\;\;
\;\;\;\;\;\;\;\;\;\;\;\;\;\;\;\;\; \mbox{if}\;\;y=\infty\; \mbox{and}\;x\in X,\\
\displaystyle\;\;\;\;\;\;\;\;\;\;0,
\;\;\;\;\;\;\;\;\;\;\;\;\;\;\;\;\;\;\;\;\;\;\;\;\;\;\;\;\;\;\;\mbox{if}\;\; x=\infty=y.
\end{cases}
\ee

Similarly, given a bounded quasi-metric space $(X,\rho)$ and a base point $c\in X$, we consider the space $X^c=X\setminus \{c\}$ and define the density function $\rho^c:X^c\times X^c\to[0,\infty)$ as follows
\be\label{q-2}\rho^c(x,y)=\rho^c(y,x)= \frac{\rho(x,y)}{\rho(x,c)\rho(y,c)}.\ee

Following \cite{BuSc}, given a quadruple $Q$ of four distinct points $a,b,c,d$ in a quasimetric space $(X,\rho),$ we denote the triple
$$M=(\rho(a,b)\rho(c,d),\rho(a,c)\rho(b,d),\rho(a,d)\rho(b,c))$$
as the {\it cross ratio triple} of $Q$.
Given a positive real number triple $M=(a,b,c)$, we call $M$ is a {\it multiplicative $K$-triple}, where $K\ge 1$, if the two largest members of $M$, say $a$ and $b$, coincide up to a multiplicative error of $K$, namely, $\frac{1}{K}\le\frac{a}{b}\le K.$

Next, we also need the following auxiliary result.

\begin{Lem}\label{z0}$($\cite[Lemma 5.1.2]{BuSc}$)$
Assume $\rho$ is a $K$-quasi-metric on $X$, $K\geq 1$. Then for every quadruple $Q$ of distinct points of $X$, the cross-ratio triple $M$ of $Q$ is a multiplicative $K^2$-triple.
\end{Lem}

Then we shall show that $\rho_a$ and $\rho^c$ defined as above for quasimetric spaces are also quasimetric.

\begin{lem}\label{z1} Let $(X,\rho)$ be a $K$-quasi-metric space.
\begin{enumerate}
  \item\label{q-0.0} If $X$ is unbounded with $a\in X$, then $(\dot{X},\rho_a)$ is a $K'$-quasi-metric space with $K'=4K^2$;
  \item\label{q-0.1} If $X$ is bounded with $c\in X$, then $(X^c,\rho^c)$ is a $K''$-quasi-metric space with $K''=K^2$.
\end{enumerate}
\end{lem}
\bpf Our proofs is a mimic of \cite[Proposition 5.3.6]{BuSc}. For completeness we show the details. To prove (\ref{q-0.0}), we extend $\rho$ on $\dot{X}\times \dot{X}$ as follows:
$$\widehat{\rho}(\infty,\infty)=0,\;\;\widehat{\rho}(x,\infty)=1+\rho(x,a)=\widehat{\rho}(\infty,x)\;\;\mbox{and}\;\;\;\;\widehat{\rho}(x,y)=\rho(x,y)$$
for all $x,y\in X$. Then we claim that $(\dot{X},\widehat{\rho})$ is a $2K$-quasi-metric space, that is,
\be\label{q-0.2}\widehat{\rho}(x,y)\leq 2K(\widehat{\rho}(x,z)\vee \widehat{\rho}(z,y))\ee
for all $x,y,z\in \dot{X}$. To this end, we consider three possibilities. If $\infty\not\in\{x,y,z\}$, then (\ref{q-0.2}) follows since $(X,\rho)$ is a $K$-quasi-metric space. For the second possibilities, $z=\infty$,  we observe that
\beq\nonumber \widehat{\rho}(x,y)&=& \rho(x,y)\leq K(\rho(x,a)\vee \rho(a,y))
\\\nonumber&\leq&   K[(1+\rho(x,a))\vee(1+\rho(a,y))]
\\\nonumber&=& K[\widehat{\rho}(x,z)\vee \widehat{\rho}(z,y)]
,\eeq
as desired. For the remaining possibility, by symmetry,  we may assume  $y=\infty$.  A direct computation gives
\beq\nonumber \widehat{\rho}(x,y)&=& 1+\rho(x,a)\leq 1+K(\rho(x,z)\vee\rho(z,a))
\\\nonumber&\leq&   1+K\rho(x,z)+K\rho(z,a)
\\\nonumber&\leq& K\widehat{\rho}(x,z)+K \widehat{\rho}(z,y)
\\\nonumber&\leq& 2K[\widehat{\rho}(x,z)\vee \widehat{\rho}(z,y)]
,\eeq
which deduces (\ref{q-0.2}).

Next for any points $x, y, z$ in $(\dot{X},\rho_a)$ be given. If one of them is $\infty$, then the required assertion for (\ref{q-0.0}) is easy to verify and we only need to  consider the case that none of them is $\infty$. It follows from Lemma \Ref{z0} that for points $x, y, z$ and $\infty$ in $(\dot{X},\widehat{\rho})$, the cross-ratio triple $M(x,y,z,\infty)$ is a multiplicative $4K^2$-triple and thus  we have
$$\widehat{\rho}(x,y)\widehat{\rho}(z,\infty)\leq 4K^2[\widehat{\rho}(x,z)\widehat{\rho}(y,\infty)\vee\widehat{\rho}(y,z)\widehat{\rho}(x,\infty)].$$
This yields
$$\rho(x,y)(1+\rho(z,a))\leq 4K^2[\rho(x,z)(1+\rho(y,a))\vee \rho(y,z)(1+\rho(x,a))]$$
and so we obtain
$$\rho_a(x,y)\leq 4K^2[\rho_a(x,z)\vee \rho_a(z,y)],$$
this proves (\ref{q-0.0}).

It remains to show (\ref{q-0.1}). Fix $x, y, z\in X^{c}$. Since $(X,\rho)$ is a $K$-quasi-metric space, again appealing to Lemma \Ref{z0}, we know that the cross-ratio triple $M(x,y,z,c)$ for the quadruple $\{x,y,z,c\}$ is a multiplicative $K^2$-triple. Thus we have
$$\rho(x,y)\rho(z,c)\leq K^2[\rho(x,z)\rho(y,c) \vee \rho(z,y)\rho(x,c)]$$
which deduces
$$\rho^c(x,y)=\frac{\rho(x,y)}{\rho(x,c)\rho(y,c)}\leq K^2[\frac{\rho(x,z)}{\rho(x,c)\rho(z,c)}\vee \frac{\rho(z,y)}{\rho(z,c)\rho(y,c)}]=K^2[\rho^c(x,z)\vee \rho^c(z,y)].$$
Hence Lemma \ref{z1} follows.
\epf

We conclude this part with the following lemma for later use.

\begin{lem}\label{z-10}
Let $(Z,\rho)$ be a bounded quasi-metric space with $c\in Z$ and $\diam (Z,\rho)=T$. Then the identity map $\varphi:(Z^c,\rho)\to (Z^c,(\rho^c)_\infty)$ is $(1+T)^2$-bilipschitz.
\end{lem}
\bpf We first flatten the quasi-metric space $(Z,\rho)$ with respect to $c\in Z$ and thus obtain the deformed space $(Z^c,\rho^c)$. Consider the one-point extension space $Z^c\cup\{\infty\}=\dot{Z^c}$. Then we extend the density function $\rho^c$ to $\dot{Z^c}$, given by
$$\rho^c(x,\infty)=\frac{1}{\rho(x,c)}\;\;\;\; \mbox{and}\;\;\;\; \rho^c(\infty,\infty)=0.$$
Note that $\rho^c$ is still a quasi-metric in $\dot{Z^c}$.

Next, we consider the sphericalization space $(\dot{Z^c},(\rho^c)_\infty)$ with respect to $\infty\in \dot{Z^c}$. For all $x,y\in Z^c$, define
$$\rho^c(x,y)=\frac{\rho(x,y)}{\rho(x,c)\rho(y,c)}\;\;\;\;\mbox{and}\;\;\;\;\rho^c(x,\infty)=\frac{1}{\rho(x,c)}.$$
Thus we have
$$\rho^c_\infty(x,y)=\frac{\rho^c(x,y)}{[1+\rho^c(x,\infty)][1+\rho^c(y,\infty)]}=\frac{\rho(x,y)}{[1+\rho(x,c)][1+\rho(y,c)]}.$$
Hence this shows Lemma \ref{z-10} because $\diam (Z,\rho)=T$.
\epf

\subsection{Sphericalization and flattening of measures}
For a quasi-metric space $(X,\rho)$ which admits a Borel regular measure $\mu$ with dense support, we form two new transforming measures under sphericalization and flattening, respectively. Let $(X,\rho,\mu)$ be a quasi-metric measure space.
If $X$ is unbounded with $a\in X$, then the spherical measure $\mu_a$ is given by
\be\label{q-1.0}\mu_a(A)=\int_{A\setminus \{\infty\}}\frac{1}{\mu(B(a,1+\rho(a,z)))^2}d\mu(z).\ee
If $X$ is bounded, then the corresponding flattening measure $\mu^c$ is defined by
\be\label{q-1.1}\mu^c(A)=\int_{A} \frac{1}{\mu(B(c,\rho(c,z)))^2}d\mu(z),\ee
where $A\subset X^c$ is a Borel set.

We next discuss the  Ahlfors regularity and doubling  properties of quasi-metric measure spaces under sphericalization and flattening transformation. In particular, we generalize \cite[Propositions 3.1 and 4.1]{LS} to quasimetric measure spaces, see Theorems \ref{z-1} and \ref{z-2} below.

\begin{thm}\label{z-1} Let $Q>0$ and $K\geq 1$. Suppose that $(X,\rho,\mu)$ is an unbounded Ahlfors $Q$-regular $K$-quasi-metric space with $a\in X$, then the sphericalization of $X$ (with respect to $a$) $(\dot{X},\rho_a,\mu_a)$ is a bounded Ahlfors $Q$-regular $K'$-quasi-metric measure space.
\end{thm}

\bpf By the definition of $\rho_a$ we see that $\diam_a \dot{X}=\sup \{\rho_a(x,y):x,y\in \dot{X}\}\leq K$. This together with Lemma \ref{z1} assert that $(\dot{X},\rho_a)$ is a bounded $K'$-quasi-metric space with $K'=4K^2$. So we only need to verify the Ahlfors $Q$-regularity of $(\dot{X},\rho_a,\mu_a)$.

Since $(X,\rho,\mu)$ is Ahlfors $Q$-regular, for all $z\in X$ we have
\be\label{q-1.2}\frac{1}{C_A}(1+\rho(z,a))^Q\leq \mu(B(a,1+\rho(a,z)))\leq C_A(1+\rho(z,a))^Q,\ee
where $C_A$ is the Ahlfors regularity constant. As the Ahlfors regularity can be  preserved under bilipschitz transformations, we see from  (\ref{q-1.0}) and (\ref{q-1.2}) that we may define the spherical measure $\mu_a$ as follows:
\be\label{q-1.3}\mu_a(B_a(x,r))=\int_{B_a(x,r)\setminus\{\infty\}}\frac{d\mu(y)}{[1+\rho(y,a)]^{2Q}}\ee
for all $x\in \dot{X}$ and $0<r\leq \diam_a \dot{X}\leq K$, where $B_a(x,r):=\{y\in\dot{X}:\rho_a(x,y)<r\}$. It follows from the definition of $\rho_a$ that $B_a(x,r)\setminus\{\infty\}$ is a Borel set of $(X,\rho)$ because the identity map $(X,\rho)\to (X,\rho_a)$ is locally bilipschitz and so the topology induced by these two quasi-metrics are coincide.

Next, we claim that 
\be\label{q-1.4}\mu_a(B_a(x,K))=\mu_a(\dot{X})\leq 4^{Q}C_A<\infty\ee
for all $x\in X$. Indeed, since $\diam_a \dot{X}\leq K$, we know that for all $x\in X$,
\beq\nonumber
\mu_a(B_a(x,K))&=&\mu_a(\dot{X}) =\int_X \frac{d\mu(y)}{[1+\rho(y,a)]^{2Q}}
\\ \nonumber &=& \left(\sum_{i=1}^\infty \int_{B(a,2^i)\setminus B(a,2^{i-1})}  +\int_{B(a,1)}\right)\frac{d\mu(y)}{[1+\rho(y,a)]^{2Q}}
\\ \nonumber &\leq& \sum_{i=1}^\infty \frac{\mu(B(a,2^i))}{(2^{i-1})^{2Q}}\leq \sum_{i=1}^\infty C_A\frac{2^{iQ}}{2^{2Q(i-1)}}
\\ \nonumber &\leq&  \frac{4^Q}{2^Q-1} C_A<\infty,\eeq
as desired. Thus we only need to find a constant $C=C(C_A,Q,K)$ such that for all $x\in X$ and $0<r\leq \frac{1}{2}$,
\be\label{zz0}\frac{1}{C}r^Q\leq \mu_a(B_a(x,r))\leq Cr^Q.\ee
Since for $r>1/2$, $\mu_a(B_a(x,r))$ must be compariable with a constant by   (\ref{zz0}) and (\ref{q-1.4}).

In the following, we split into several cases to prove (\ref{zz0}). We begin this discussion with a consideration of the balls centered at $\infty$.
\bca\label{zz1} $x=\infty$.   \eca
Set $R=\frac{1}{r}-1$. Thus we have $\frac{1}{2r}\leq R<\frac{1}{r}$ because $r\leq \frac{1}{2}$. Note that
\beq\nonumber
B_a(\infty,r)\setminus\{\infty\}&=& \{y\in X:\rho_a(y,\infty)=\frac{1}{1+\rho(x,a)}<r\}
\\ \nonumber &=& \{y\in X:\rho(y,a)> R\}=X\setminus \overline{B}(a,R).
\eeq
Let $\kappa=2C_A^{2/Q}$ and $B=B(a,\kappa^iR)$ ($i=0,1,2,...$). Then one computes
$$\mu(B_i\setminus B_{i-1})=\mu(B_i)-\mu(B_{i-1})\geq \frac{(\kappa^iR)^Q}{C_A}-C_A(\kappa^{i-1}R)^Q\geq \frac{(\kappa^iR)^Q}{2C_A},$$
and so
\beq\nonumber
\mu_a(B_a(\infty,r)) &=& \int_{X\setminus \overline{B}(a,R)} \frac{d\mu(y)}{[1+\rho(a,y)]^{2Q}}= \sum_{i=1}^\infty \int_{B_i\setminus B_{i-1}} \frac{d\mu(y)}{[1+\rho(a,y)]^{2Q}}
\\ \nonumber &\leq& \sum_{i=1}^\infty \frac{\mu(B_i)}{(\kappa^{i-1}R)^{2Q}}\leq \sum_{i=1}^\infty \frac{(\kappa^iR)^Q}{(\kappa^{i-1}R)^{2Q}}
\\ \nonumber &\leq& {\frac{2C_A\kappa^{2Q}}{R^Q}<2^{Q+1}C_A\kappa^{2Q}r^Q,}
\eeq
On the other hand, we get
$$\mu_a(B_a(\infty,r))\geq \sum_{i=1}^\infty \frac{\mu(B_i\setminus B_{i-1})}{(2\kappa^{i}R)^{2Q}}\geq  \sum_{i=1}^\infty \frac{R^Q\kappa^{iQ}}{2C_A(2\kappa^{i}R)^{2Q}}\geq \frac{1}{2C_A(4\kappa)^{Q}R^Q}\geq \frac{r^Q}{2C_A(4\kappa)^{Q}},$$
these estimates imply (\ref{zz0}).

\bca\label{zz2} $\rho_a(x,\infty)\leq \frac{r}{K'}$.   \eca
In this case, we claim that
\be\label{q-1.5} B_a(\infty,\frac{r}{K'})\subset B_a(x,r)\subset B_a(\infty,K'r).\ee
This can be seen as follows. For all $y\in B_a(\infty,\frac{r}{K'})$, we have $\rho_a(x,y)\leq K'(\rho_a(x,\infty)\vee \rho_a(\infty,y))<r$. On the other hand, for all $z\in B_a(x,r)$,  we get $\rho_a(\infty,z)\leq K'(\rho_a(x,\infty)\vee \rho_a(x,z))<K'r$, as claimed.

Consequently, combining the inclusion relation (\ref{q-1.5}) with Case \ref{zz1}, we obtain (\ref{zz0}).

\bca $\frac{r}{K'}<\rho_a(x,\infty)\leq 4K^2r$.   \eca
A similar argument as Case \ref{zz2}, we can obtain the inclusion $B(x,\frac{1}{16K^4K'r})\subset B_a(x,r)\subset B_a(\infty,4K'K^2r)$. Indeed, for all $y\in B(x,\frac{1}{16K^4K'r})$, we have
\beq\nonumber
\rho_a(x,y) &=&\rho(x,y)\rho_a(x,\infty)\rho_a(y,\infty)<\frac{4K^2r}{16K^4K'r}\rho_a(y,\infty)
\\ \nonumber &\leq& \frac{K'}{4K^2K'}(\rho_a(x,y)\vee \rho_a(x,\infty))\leq \frac{1}{4K^2}(\rho_a(x,y)\vee 4K^2r)<r,
\eeq
and for all $z\in B_a(x,r)$,
$$\rho_a(z,\infty)\leq K'(\rho_a(z,x)\vee \rho_a(x,\infty))<K'(r\vee 4K^2r)=4K'K^2r,$$
as required. Moreover, the desired upper Ahlfors $Q$-regularity can be obtained by means of the inclusion $B_a(x,r)\subset B_a(\infty,4K'K^2r)$ and the implication in Case \ref{zz1}. To show the lower regularity, for all  $y\in B(x,\frac{1}{16K^4K'r})$ we compute
$$\rho_a(y,\infty)=\frac{\rho_a(x,y)}{\rho(x,y)\rho_a(x,\infty)}>4K^2K'\rho_a(x,y)>\rho_a(x,y)$$
and so $\frac{r}{K'}\leq \rho_a(x,\infty)\leq K'(\rho_a(x,y)\vee \rho_a(y,\infty))=K'\rho_a(y,\infty)$. This yields
\beq\nonumber
\mu_a(B_a(x,r))&=&\int_{B_a(x,r)\setminus \{\infty\}} \frac{d\mu(y)}{[1+\rho(y,a)]^{2Q}}
\\ \nonumber &\geq& \int_{B(x,\frac{1}{16K^4K'r})} \frac{d\mu(y)}{[1+\rho(y,a)]^{2Q}}
\\ \nonumber &\geq& \mu(B(x,\frac{1}{16K^4K'r}))(\frac{r}{K^{'2}})^{2Q}
\\ \nonumber &\geq& \frac{1}{C_A}(\frac{1}{16K^4K'r})^Q(\frac{r}{K^{'2}})^{2Q}=\frac{1}{C}r^Q.
\eeq
Therefore, we are done in this case.

\bca\label{zz3} $\rho_a(x,\infty)\geq 4K^2r$ and $\rho(x,a)\leq 1$. \eca
In this case, we are going to build the inclusion
\be\label{q-1.6} B(x,r)\subset B_a(x,r)\subset B(x,8Kr).\ee
Towards this end, for all $y\in B(x,r)$ we have
$$\rho_a(x,y)= \rho(x,y)\rho_a(x,\infty)\rho_a(y,\infty)\leq \rho(x,y)<r$$
and for all $z\in B_a(x,r)$,
$$r>\rho_a(x,z)= \rho(x,z)\rho_a(x,\infty)\rho_a(z,\infty)\geq 4K^2r\rho(x,z)\rho_a(z,\infty).$$
From which we deduce that
$$4K^2\rho(x,z)\leq\frac{1}{\rho_a(z,\infty)} = 1+\rho(a,z) \stackrel{(*)}{\leq} 2K[\rho(x,z)\vee (1+\rho(x,a))]=2K(1+\rho(x,a))<4K, $$
where the inequality $(*)$ follows from a similar argument as  in Lemma \ref{z1}. Thus we obtain $1+\rho(a,z)\leq 4K$ and $1+\rho(a,x)\leq 2$. Moreover, we get
$$\rho(x,z)=\rho_a(x,z)[1+\rho(x,a)][1+\rho(z,a)]<8Kr,$$
this yields (\ref{q-1.6}).

Then by (\ref{q-1.6}) we compute
$$\mu_a(B_a(x,r))=\int_{B_a(x,r)}\frac{d\mu(y)}{[1+\rho(y,a)]^{2Q}}\geq \frac{\mu(B_a(x,r))}{(4K)^{2Q}}\geq \frac{\mu(B(x,r))}{(4K)^{2Q}}\geq \frac{r^Q}{C_A(4K)^{2Q}}.$$
Moreover, again by (\ref{q-1.6}) we get
$$\mu_a(B_a(x,r))\leq \mu(B_a(x,r))\leq \mu(B(x,8Kr))\leq C_A(8Kr)^Q,$$
and hence we prove (\ref{zz0}) by means of these two estimates.

\bca $\rho_a(x,\infty)\geq 4K^2r$  and $\rho(x,a)> 1$. \eca
Put $t_y=1+\rho(a,y)$ for $y\in B_a(x,r)$. Thus we claim that
\be\label{q-1.7} \frac{t_x}{2K}<t_y<2Kt_x\;\;\;\;\;\;\mbox{and}\;\;\;\;\;\; B(x,\frac{r}{2K}t_x^2)\subset B_a(x,r)\subset B(x,2Krt_x^2).\ee
Indeed, since $\rho_a(x,\infty)\geq 4Kr$, a similar argument as Case \ref{zz3} gives $4K^2\rho(x,y)<1+\rho(a,y)$ for all $y\in B_a(x,r)$. And since $(\rho(x,y),1+\rho(x,a),1+\rho(y,a))$ is a $2K$-triple, we know that
$$4K^2\rho(x,y)<1+\rho(y,a)\leq 2K[\rho(x,y)\vee (1+\rho(x,a))]=2K(1+\rho(x,a))$$
and so $t_y<2Kt_x$ and $2K\rho(x,y)\leq 1+\rho(x,a)$. On the other hand, since
$$t_x=1+\rho(x,a)\leq 2K[\rho(x,y)\vee (1+\rho(y,a))]=2K(1+\rho(y,a))=2Kt_y,$$
we deduce the first inequalities $\frac{t_x}{2K}<t_y<2Kt_x$ of (\ref{q-1.7}).

It remains to show the inclusions of (\ref{q-1.7}). For all $y\in B_a(x,r)$, we have
$$\rho(x,y)=\rho_a(x,y)t_xt_y<2Kt_x^2r,$$
which implies the right hand side inclusion. For each $z\in B(x,\frac{r}{2K}t_x^2)$, we have
$$t_x\leq 2K[\rho(x,z)\vee (1+\rho(z,a))]\leq 2K[\frac{r}{2K}t_x^2 \vee (1+\rho(z,a))]=2K(1+\rho(a,z)),$$
because $\rho_a(x,\infty)\geq 4K^2r$ implies $rt_x\leq \frac{1}{4K^2}$. Moreover, we obtain
$$\rho_a(x,z)=\frac{\rho(x,z)}{t_x(1+\rho(a,z))}<\frac{\frac{r}{2K}t_x^2}{t_x\cdot \frac{t_x}{2K}}=r,$$
as required. Consequently, we get (\ref{q-1.7}).

Now, we shall appeal (\ref{q-1.7}) and the Ahlfors $Q$-regularity of $(X,\rho,\mu)$ to obtain
\beq\nonumber
\mu_a(B_a(x,r))&=&\int_{B_a(x,r)\setminus \{\infty\}} \frac{d\mu(y)}{[1+\rho(y,a)]^{2Q}}<(\frac{2K}{t_x})^{2Q}\mu(B_a(x,r))
\\ \nonumber &\leq& (\frac{2K}{t_x})^{2Q}\mu(B(x,2Krt_x^2))
\\ \nonumber &<& (\frac{2K}{t_x})^{2Q}C_A(2Krt_x^2)^Q=Cr^Q
\eeq
and
$$\mu_a(B_a(x,r))>\frac{\mu(B_a(x,r))}{(2Kt_x)^{2Q}}>\frac{\mu(B(x,\frac{r}{2K}t_x^2))}{(2Kt_x)^{2Q}}>\frac{(\frac{r}{2K}t_x^2)^Q}
{C_A(2Kt_x)^{2Q}}>\frac{r^Q}{C}.$$

Hence we complete the proof of Theorem \ref{z-1}.
\epf

\begin{thm}\label{z-2} Let $Q>0$ and $K\geq 1$. Suppose that $(X,\rho,\mu)$ is a bounded Ahlfors $Q$-regular $K$-quasi-metric space with $c\in X$, then the flattening of $X$ with respect to $c$, $(X^c,\rho^c,\mu^c)$ is an unbounded Ahlfors $Q$-regular $K''$-quasi-metric measure space.
\end{thm}

\bpf First, from Lemma \ref{z1} it follows that $(X^c,\rho^c)$ is a $K''$-quasi-metric space with $K''=K^2$. Thus it suffices to show that $(X^c,\rho^c,\mu^c)$ is Ahlfors $Q$-regular. To this end, since $(X,\rho,\mu)$ is Ahlfors $Q$-regular, for all $z\in X^c=X\setminus \{c\}$ we have
$$\frac{1}{C_A}\rho(c,z)^Q\leq \mu(B(c,\rho(c,z)))\leq C_A\rho(c,z)^Q,$$
where $C_A$ is the Ahlfors regularity constant. Because bilipschitz homeomorphism preserves Ahlfors $Q$-regularity, we may assume that the flattening measure $\mu^c$ on $(X^c,\rho^c)$ is given by
$$\mu^c(A)=\int_A\frac{d\mu(z)}{\rho(c,z)^{2Q}},$$
where $A\subset X^c$. Fix $x\in X^c$ and $r>0$, it suffices to show that $\mu^c(B^c(x,r))$ is compariable with $r^Q$, where $B^c(x,r)=\{y\in X^c:\rho^c(x,y)<r\}$. We consider three cases.

{\bf Case A:} $r\rho(x,c)\leq \frac{1}{2K}$.
In this case, we first claim that
\be\label{q-2.0} B(x,\frac{r}{K}\rho(x,c)^2)\subset B^c(x,r)\subset B(x,Kr\rho(x,c)^2).\ee
Indeed, for each $y\in B^c(x,r)$, we get
$$\rho(y,x)=\rho^c(x,y)\rho(x,c)\rho(y,c)\leq \frac{1}{2K}\rho(y,c)\leq \frac{1}{2}(\rho(y,x)\vee \rho(x,c))=\frac{1}{2}\rho(x,c),$$
which deduces $\rho(x,y)\leq Kr\rho(x,c)^2$. Moreover, for every $z\in B(x,\frac{r}{K}\rho(x,c)^2)$, we have
\beq\nonumber
\rho(x,c) &\leq& K(\rho(x,z)\vee \rho(z,c))\leq K(\frac{r}{K}\rho(x,c)^2\vee \rho(z,c))
\\ \nonumber&\leq& K(\frac{1}{2K^2}\rho(x,c)\vee \rho(z,c))=K\rho(z,c),
\eeq
and so $\rho^c(x,z)=\frac{\rho(x,z)}{\rho(x,c)\rho(c,z)}<\frac{\frac{r}{K}\rho(x,c)}{\rho(c,z)}<r$, which implies (\ref{q-2.0}).

Next, we shall use (\ref{q-2.0}) to estimate $\mu^c(B^c(x,r))$. On one hand, one computes
\beq\nonumber
\mu^c(B^c(x,r))&=&\int_{B^c(x,r)} \frac{d\mu(y)}{\rho(c,y)^{2Q}}\geq \frac{\mu(B^c(x,r))}{[K\rho(x,c)]^{2Q}}
\\ \nonumber &\geq& \frac{\mu(B(x,\frac{r}{K}\rho(x,c)^2))}{[K\rho(x,c)]^{2Q}}
\\ \nonumber &\geq& \frac{1}{C_A}\frac{[\frac{r}{K}\rho(x,c)^2)]^Q}{[K\rho(x,c)]^{2Q}}=\frac{r^Q}{C_AK^{3Q}}.
\eeq
On the other hand, we obtain
$$\mu^c(B^c(x,r))\leq \frac{\mu(B^c(x,r))}{(\rho(x,c)/K)^{2Q}}\leq \frac{\mu(B(x,Kr\rho(x,c)^2))}{(\rho(x,c)/K)^{2Q}}\leq C_AK^{3Q}r^Q.$$
Hence we are done in this case.

{\bf Case B:} $r\rho(x,c)\geq 2KC_A^2$.
We first establish the inclusions which is needed in the estimation of the measure $\mu^c(B^c(x,r))$, that is,
\be\label{q-2.1} X\setminus \overline{B}(c,\frac{K}{r})\subset B^c(x,r)\subset X\setminus \overline{B}(c,\frac{1}{K r}).\ee
To this end, for all $y\in X^c$ with $\rho(y,c)>K/r$, we compute
$$\rho^c(x,y)=\frac{\rho(x,y)}{\rho(x,c)\rho(c,y)}\leq \frac{K(\rho(x,c)\vee \rho(c,y))}{\rho(x,c)\rho(c,y)}=\frac{K}{\rho(x,c)\wedge \rho(c,y)}<\frac{K}{\frac{2KC_A^2}{r}\wedge \frac{K}{r}}=r.$$
Moreover, for every $z\in B^c(x,r)$, $\rho(x,z)=\rho^c(x,z)\rho(x,c)\rho(z,c)<r\rho(x,c)\rho(z,c)$, and so
$$\rho(x,c)\leq K(\rho(x,z)\vee \rho(z,c))\leq K(r\rho(x,c)\vee 1)\rho(z,c)=Kr\rho(x,c)\rho(z,c),$$
which implies $\rho(z,c)\geq \frac{1}{Kr}$. Consequently, from the above estimates (\ref{q-2.1}) follows.

Next, we are going to estimate $\mu^c(B^c(x,r))$ by virtue of (\ref{q-2.1}). For the upper Ahlfors regularity,  appealing to the regularity of the space $(X,\rho,\mu)$, we get
\beq\nonumber
\mu^c(B^c(x,r)) &=& \int_{B^c(x,r)} \frac{d\mu(z)}{\rho(z,c)^{2Q}}\leq \int_{X\setminus \overline{B}(c,\frac{1}{K r})} \frac{d\mu(z)}{\rho(z,c)^{2Q}}
\\ \nonumber&=& \sum_{i=0}^\infty \int_{B(c,\frac{2^{i+1}}{K r})\setminus B(c,\frac{2^{i}}{K r})} \frac{d\mu(z)}{\rho(z,c)^{2Q}}
\\ \nonumber&\leq& \sum_{i=0}^\infty\frac{\mu(B(c,\frac{2^{i+1}}{K r}))}{(\frac{2^{i}}{K r})^{2Q}}\leq \sum_{i=0}^\infty \frac{C_A(\frac{2^{i+1}}{K r})^Q}{(\frac{2^{i}}{K r})^{2Q}}
\\ \nonumber&\leq& C_A(2K)^{Q+1} r^Q.
\eeq
It remains to verify the lower regularity. A direct computation gives
\beq\nonumber
\mu^c(B^c(x,r))&=&\int_{B^c(x,r)} \frac{d\mu(y)}{\rho(c,y)^{2Q}}\geq \int_{X\setminus \overline{B}(c,\frac{K}{r})} \frac{d\mu(z)}{\rho(z,c)^{2Q}}
\\ \nonumber &\geq&  \int_{B(c,\frac{2KC_A^2}{r})\setminus \overline{B}(c,\frac{K}{r})} \frac{d\mu(z)}{\rho(z,c)^{2Q}}
\\ \nonumber &\geq& \frac{\mu\big(B(c,\frac{2KC_A^2}{r})\setminus \overline{B}(c,\frac{K}{r})\big)}{(\frac{2KC_A^2}{r})^{2Q}}
\\ \nonumber &\geq& \frac{ \frac{1}{C_A}(\frac{2KC_A^2}{r})^Q- C_A(\frac{K}{r})^Q}{(\frac{2KC_A^2}{r})^{2Q}}>\frac{r^Q}{C},
\eeq
as desired.

{\bf Case C:} $\frac{1}{2K}\leq r\rho(x,c)\leq 2KC_A^2$.
In this case, we see that $B^c(x,\frac{r}{4K^2C_A^2})$ satisfies the assumption of Case $A$ and similarly $B^c(x, 4K^2C_A^2r)$ satisfies the condition of Case $B$. Both of them, their measures with respect to the flattening measures, are compariable with $r^Q$. Hence the same assertion holds also for $B^c(x,r)$.

Hence we prove Theorem \ref{z-2}.
\epf
%
\section{Preservation of doubling measure}
Recently, the third author and Shanmugalingam proved that both the metric space flattening and sphericalization  preserving doubling measure, see \cite[Propositions 3.3 and 4.2]{LS}. In this section, we mainly generalize these results to quasi-metric measure spaces and show that doubling property for uniformly perfect quasi-metric measure spaces is preserved under flattening and sphericalization transformations. Our proofs are based on Theorems \ref{z-1} and \ref{z-2} and the generalized David-Semmes's deformation\cite{DS}. Our main result in this section  is stated as follows.
\begin{thm}\label{z-6}
Let $(X,\rho)$ be a $\tau$-uniformly perfect $K$-quasi-metric space which carries a doubling measure $\mu$ with coefficient $C_\mu \geq 1$.
\begin{enumerate}
  \item If $X$ is bounded with $c\in X$, then the flattening measure $\mu^c$ on $(X^c,\rho^c)$ is doubling;
  \item If $X$ is unbounded with $a\in X$, then the spherical measure $\mu_a$ on $(\dot{X},\rho_a)$ is doubling.
\end{enumerate}
\end{thm}

\subsection{Auxiliary results} In this subsection, we are going to prove some auxiliary results for later use.
\begin{lem}\label{z-5}
Let $(Z,\rho)$ be a $\tau$-uniformly perfect $K$-quasi-metric space which carries a doubling measure $\mu$ with coefficient $C\geq 1$. Then there are constants $\alpha>0$ and $C_0>0$ such that
$$\frac{\mu(B(x,r))}{\mu(B(x,R))}\leq C_0(\frac{r}{R})^\alpha,$$
for all $x\in Z$ and $0<r\leq R\leq \diam Z$.
\end{lem}

\bpf We first show that there are constants $\delta_1,\delta_2\in(0,1)$ such that
\be\label{ww0} \mu(B(x,\delta_1r))\leq \delta_2\mu(B(x,r))\ee
for all $x\in Z$ and $0<r \leq \diam(Z)$. Since $Z$ is $\tau$-uniformly perfect, there is a point $y\in Z$ such that
$$\frac{\tau r}{4K}\leq \rho(x,y)\leq \frac{r}{2K}.$$
Thus we know that $B(x,\frac{\tau r}{8K^2})$ and $B(y,\frac{\tau r}{8K^2})$ are disjoint sets of $B(x,r)$, which implies
$$\mu(B(x,\frac{\tau r}{8K^2}))\leq \mu(B(x,r))- \mu(B(y,\frac{\tau r}{8K^2}))\leq (1-\frac{1}{C^{\log_2\frac{K^3}{\tau}+4}})\mu(B(x,r)),$$
since $(Z,\mu)$ is a $C$-doubling measure space. This yields $(\ref{ww0})$ with the choice of $\delta_1=\frac{\tau }{8K^2}$ and $\delta_2=1-\frac{1}{C^{\log_2\frac{K^3}{\tau}+4}}$.

For $0<r\leq R\leq \diam(Z)$, there is an integer $n$ such that $\delta_1^n<r/R\leq \delta_1^{n-1}$. Let $\alpha=\log_{\delta_1}\delta_2$ and $C_0=1/\delta_2$. By $(\ref{ww0})$ we obtain
\begin{eqnarray*} \mu(B(x,r)) &\leq& \mu(B(x,\delta_1^{n-1}R))\leq \delta_2^{n-1}\mu(B(x,R))
\\ \nonumber&\leq& \delta_2^{\log_{\delta_1}\frac{r}{R}-1}\mu(B(x,R))=C_0(\frac{r}{R})^{\alpha}\mu(B(x,R)).
\end{eqnarray*}
Hence Lemma \ref{z-5} follows.
\epf

\begin{lem}\label{zz4} Let $f:(Z,\rho)\to (Z',\rho')$ be an $\eta$-quasisymmetric homeomorphism between two quasi-metric spaces. Then for any $k\geq 1$, $x\in Z$ and $0<r\leq \diam(Z)$, there exists a ball $B(x',R)$ centered at $x'=f(x)$ with radius $R$ such that
$$B(x',R)\subset f(B(x,r)) \subset  f(B(x,kr)) \subset B(x',\eta(k)R).$$
\end{lem}

\bpf Since $0<r\leq \diam Z$, we may assume that $Z\setminus B(x,r)\neq \emptyset$. Let
$$R=\inf\{\rho'(x',f(z)): \; z\in Z\setminus B(x,r)\}.$$
We claim that $$B(x',R)\subset f(B(x,r)).$$
Otherwise, there is a point $x_0'\in B(x',R)$ with $x_0'=f(x_0)$ and $x_0\in Z\setminus B(x,r)$. Then by the choice of $R$, we have $\rho'(x_0',x')<R\leq \rho'(x_0',x')$,  which conducts a contradiction. Next, we show that for all $k\ge 1$
$$f(B(x,kr)) \subset B(x',\eta(k)R).$$

Note that for all $z\in B(x,kr)$ and $w\in Z\setminus B(x,r)$, we have
$\rho(x,z)\leq kr\leq k\rho(x,w)$ and so $$\rho'(x',f(z))\leq \eta(k)\rho(x',f(w)),$$
since $f$ is $\eta$-quasisymmetric. From which  we obtain above, the Lemma \ref{zz4} follows.
\epf

\begin{lem}\label{z-9}
 Suppose that a quasi-metric space is quasi-symmetrically embedding into a quasi-metric space which carries a doubling measure, then the pull-back measure on the pre-image space is also doubling.
\end{lem}
\bpf Assume that $f:(Z,\rho)\to (Z',\rho')$ is an $\eta$-quasi-symmetric embedding between two quasi-metric spaces and $\nu$ is a doubling measure on $(Z',\rho')$. Define the pull-back measure induced by $f$
$$\mu_f(E):=\nu(E)$$
for any Borel set $E\subset Z$. We need to show that $\mu_f$ on $(Z,\rho)$ is  doubling Borel measure.

To this end, we observe that every quasi-symmetric mapping is a Borel function and a similar argument as \cite[$3.3.21$]{HKST} shows that $\mu_f$ is a Borel measure. It remains to verify the doubling property of $\mu_f$. For every $z\in Z$ and $0<r\leq \diam Z$, it follows from Lemma \ref{zz4} that there is a positive number $R$ such that
$$B(f(x),R)\subset f(B(x,r)) \subset  f(B(x,2r)) \subset B(f(x),\eta(2)R).$$
Therefore, we compute
\beq\nonumber
\mu_f(B(x,2r)) &=& \nu(f(B(x,2r)))\leq \nu(B(f(x),\eta(2)R))
\\ \nonumber&\leq& C_\nu^{\log_2\eta(2)+1}\nu(B(f(x),R))
\\ \nonumber&\leq& C_\nu^{\log_2\eta(2)+1}\mu_f(B(x,r)).
\eeq
This proves Lemma \ref{z-9}.
\epf

Next, we generalize  the David-Semmes's deformation for doubling quasi-metric measure spaces.
\begin{thm}\label{z-3}
A quasi-metric space is quasisymmetrically equivalent to an Ahlfors regular space if and only if it is an uniformly perfect doubling measure space.
\end{thm}
We remark our proof is a minic of \cite{S} or \cite{DS}. For completeness, we show the details in the following subsection.
\subsection{The proof of Theorem \ref{z-3}}
We note that every Ahlfors regular space is doubling, so the necessity easily follows from Lemma \ref{z-9} and \cite[Lemma 2.10]{WZ}. Hence, we only need to consider the sufficiency. To this end, we divide the proof into several steps. First,  we  deform the quasi-metric space via doubling measure and show that the deforming space is a quasi-metric space.
Next we shall show that the identity map between the deforming space and the origin space is quasisymmetric. Finally, with the help of Lemma \ref{z-5}, we  complete the proof of this theorem by showing Ahlfors regularity of the deforming space.

Assume that $(Z,\rho,\mu)$ is a $\tau$-uniformly perfect $K$-quasi-metric space which carries a doubling measure with coefficient $C>1$. Define

$$\beta(x,y)=\begin{cases}
\displaystyle\;\;\; \;\;\;0,\;\;\;\;\;\;\;\;
\;\;\;\mbox{if}\;\;x=y,\\
\displaystyle\mu(B_{x,y})^{\varepsilon},
\;\;\;\;\;\; \mbox{if}\;\;x\neq y,\;
\end{cases}
$$
where $\varepsilon>0$ and $B_{x,y}=B(x,\rho(x,y))\cup B(y,\rho(x,y))$. We claim that $(Z,\beta)$ is a quasimetric space.

\blem\label{q-3.0} $(Z,\beta)$ is a quasimetric space. \elem
\bpf
By the definition of $\beta$, it is not difficult to see that $\beta$ is a nonnegative symmetric function with $\beta(x,x)=0$ for $x\in Z$. Since $\mu(B_{x,y})\geq \mu(B(x,\rho(x,y)))>0$ for $x\neq y\in Z$, we have $\beta(x,y)=0$ if and only if $x=y$. So we only need to show that there exists a  constant $K_0\geq 1$ such that
$$\beta(x,y)\leq K_0(\beta(x,z)\vee\beta(z,y))$$
for any $x,y,z\in Z$. Without loss of generality, we may assume that $\beta(x,z)\geq\beta(z,y)$ and $x\neq y$. A direct computation gives
$$B_{x,y}\subset B(x,K\rho(x,y))\subset B(x,K^2\rho(x,z)).$$
This, together with the assumption ``$(X,\rho,\mu)$ is $C$-doubling", yields
$$\beta(x,y)=\mu(B_{x,y})^{\varepsilon} \leq C^{(2\log_2K+2)\varepsilon}\mu(B(x,\rho(x,z)))^{\varepsilon}\leq K_0\beta(x,z), $$
where $K_0=C^{(2\log_2K+2)\varepsilon}$. Hence Lemma \ref{q-3.0} follows.
\epf

\blem\label{ww-4} The identity map $(Z,\rho)\to (Z,\beta)$ is quasisymmetric. \elem
\bpf
For any three distinct points $x,a,b\in Z$ with $\rho(x,a)\leq t \rho(x,b)$, we   divide the proof into two cases. One is $Kt\leq 1$.  With the aid of Lemma \ref{z-5}, we have
$$\frac{\beta(x,a)}{\beta(x,b)}= \left[\frac{\mu(B_{x,a})}{\mu(B_{x,b})}\right]^\varepsilon\leq \left[\frac{\mu(B(x,K\rho(x,a)))}{\mu(B(x,\rho(x,b)))}\right]^{\epsilon} \leq
\left[\frac{\mu(B(x,Kt\rho(x,b)))}{\mu(B(x,\rho(x,b)))}\right]^\varepsilon\leq C_0^{\varepsilon}(Kt)^{\alpha\varepsilon},$$
where $C_0$ and $\alpha$ are the constants in Lemma \ref{z-5}. Another case is $Kt> 1$. By a similar argument as in the first case, we get
$$\frac{\beta(x,a)}{\beta(x,b)}\leq \left[\frac{\mu(B(x,Kt\rho(x,b)))}{\mu(B(x,\rho(x,b)))}\right]^\varepsilon\leq C^{(\log_2(Kt)+1)\varepsilon},$$
the last inequality follows from the assumption ``$\mu$ is a $C$-doubling measure".

Hence Lemma \ref{ww-4} is proved.
\epf

\blem\label{q-0} $(Z,\beta,\mu)$ is Ahlfors $\frac{1}{\varepsilon}$-regular. \elem
\bpf
Denote the $\beta$-ball by
$$B_\beta(x,r)=\{y\in Z:\beta(x,y)<r\}$$
for $x\in Z$ and $0<r\leq \diam_\beta(Z)$. We first note that Lemma  \ref{ww-4} implies that the $\rho$-topology coincides with the $\beta$-topology in $Z$  and from which it follows that $B_\beta(x,r)$ is a  $\rho$-open and $\mu$-measurable subset of $Z$. Take a smallest number $s>0$ so that $B_\beta(x,r)\subset \overline{B}(x,s)$. Then there exists a point $z\in B_\beta(x,r)$ such that $\rho(x,z)>s/2$. We claim that
\be\label{ww-5} B(x,s/c_1)\subset B_\beta(x,r)\subset \overline{B}(x,s), \ee
where $c_1=2C_0C^{\log_2K+1}$ and $C_0$ is the constant in Lemma \ref{z-5}.

By the choice of $s$, we only need to verify the first inclusion. For every $y\in B(x,s/c_1)$, we have
\be\label{q-3.1} \rho(x,y)\leq \frac{s}{c_1}\leq \frac{2}{c_1}\rho(x,z).\ee
Moreover since $\mu$ is a $C$-doubling measure on $(Z,\rho)$, we have
$$\beta(x,y)=\mu(B_{x,y})^{\varepsilon}\leq \mu(B(x,K\rho(x,y)))^{\varepsilon}\leq C^{(\log_2K+1)\varepsilon}\mu(B(x,\rho(x,y)))^{\varepsilon}\leq \beta(x,z)<r,$$
the last inequality but one follows from Lemma \ref{z-5}, (\ref{q-3.1}) and the choice of $c_1$. Hence we get (\ref{ww-5}).

We continue the proof of this lemma. We see from  (\ref{ww-5}), the doubling property of $(Z,\rho,\mu)$ and the measurability of $B_\beta(x,r)$ that to get the Ahlfors regularity of the space $(Z,\beta,\mu)$, we only need to check the following double inequalities
\be\label{ww-6} c_2^{-1}r\leq \mu(B(x,s))^\varepsilon\leq c_2r, \ee
hold for some suitable constant $c_2>0$.

Since $(Z,\rho,\mu)$ is a $C$-doubling measure space, we have
$$r\geq \beta(x,z)\geq \mu(B(x,\rho(x,z)))^\varepsilon \geq \mu(B(x,s/2))^\varepsilon \geq \frac{1}{C^\varepsilon}\mu(B(x,s))^\varepsilon,$$
from which the upper bound in (\ref{ww-6}) follows. To prove the lower bound of \eqref{ww-6},  we divide the proof  into two cases.

For the first case, if $\diam(Z)\geq s/\tau$, where $\tau\in(0,1)$ is the uniformly perfect coefficient, we observe from the choice of $s$ that $\beta(x,w)\geq r$. Choose a point $w\in Z$ such that $s<\rho(x,w)\leq s/\tau$. Then it follows from the doubling property of $(Z,\rho,\mu)$ that
\beq\nonumber
r&\leq& \beta(x,w)=\mu(B_{x,w})^\varepsilon \leq \mu(B(x,K\rho(x,w)))^\varepsilon
\\ \nonumber &\leq& \mu(B(x,Ks/\tau))^\varepsilon\leq c_2\mu(B(x,s))^\varepsilon,
\eeq
where $c_2= C^{(\log_2\frac{K}{\tau}+1)\varepsilon}$.

For the remaining case, that is, $\diam(Z)< s/\tau$, we obtain from \eqref{ww-5} that
$$B(x,s/c_1)\subset B_\beta(x,r)\subset Z \subset B(x,s/\tau),$$
and so
\beq\nonumber
r&\leq& \diam_\beta(Z)= \sup\{\mu(B_{x_1,x_2})^\varepsilon :x_1,x_2\in Z\}
\\ \nonumber &\leq& \mu(Z)^\varepsilon\leq \mu(B({x,\frac{s}{\tau}}))^\varepsilon \leq c_2\mu(B({x,s}))^\varepsilon.
\eeq
This implies $(\ref{ww-6})$ and the proof of Lemma \ref{q-0} is complete.
\epf

Hence, Theorem \ref{z-3} follows from Lemmas \ref{q-3.0}, \ref{ww-4} and \ref{q-0} .
\subsection{The proof of Theorem \ref{z-6}}
$(1)$ Assume that $X$ is bounded with $c\in X$. According to Theorem \ref{z-3}, we know that there is an Ahlfors regular quasi-metric measure  space $(X,\delta,\mu)$ such that the identity map $\varphi:(X,\rho)\to (X,\delta)$ is quasi-symmetric. In the following, we consider the flattening transformation of the spaces $(X,\rho,\mu)$ and  $(X,\delta,\mu)$ with respect to the same point $c\in X$, respectively. By (\ref{q-2}), we know that the identity maps $\psi_\rho:(X^c,\rho)\to (X^c,\rho^c)$ and $\psi_\delta:(X^c,\delta)\to (X^c,\delta^c)$ are both M\"{o}bius homeomorphism with $\psi_\rho(c)=\psi_\delta(c)=\infty$. Since quasi-symmetric mapping is quasim\"{o}bius and the composition of quasim\"{o}bius maps is also quasim\"{o}bius, we find that
$$\varphi^c=\psi_\delta\circ \varphi\circ \psi_\rho^{-1}:(X^c,\rho^c)\to (X^c,\delta^c)$$
is quasim\"{o}bius with $\varphi^c(\infty)=\infty$ and so $\varphi^c$ is quasisymmetric by means of \cite[Theorem 3.20]{Vai-5}.

On the other hand, by virtue of Theorem \ref{z-2} and the Ahlfors regularity of $(X,\delta,\mu)$, we know that the flattening space $(X^c,\delta^c,\mu^c)$ is Ahlfors regular as well. Consequently, again appealing to Theorem \ref{z-3} the doubling property of $\mu^c$ on $(X^c,\rho^c)$ follows.

$(2)$ Assume that $(X,\rho)$ is unbounded and it  admits a doubling measure $\mu$. Let $a\in X$. Without loss of generality, we may normalize the situation so that $\mu(B_\rho(a,1))=1$. Then it follows from Theorem \ref{z-3} that the deformed space $(X,\delta)$ given by

$$\delta(x,y)=\begin{cases}
\displaystyle\;\;\; \;\;\;0,\;\;\;\;\;\;\;\;\;\;\;\;\;\;\;\;\;\;\;\;\;\;\;\;\;\;\;\;\;\;\;\;\;\;\;\;\;\;\;\;\;\;\;\;
\;\;\;\mbox{if}\;\;x=y,\\
\displaystyle\mu(B_\rho(x,\rho(x,y))\cup B_\rho(y,\rho(x,y))),
\;\;\;\;\;\; \mbox{if}\;\;x\neq y,\;
\end{cases}
$$
 is a $K_0$-quasi-metric space with $K_0$ depending only on $C_\mu,\tau$ and $K$, and moreover, the measure $\mu$ on $(X,\delta)$ is Ahlfors regular and the original space $(X,\rho)$ is quasi-symmetrically equivalent to $(X,\delta)$ via the identity map $\varphi$.

 Next, we    consider the sphericalization spaces of the above two spaces associated to the point $a$, that is, $(\dot{X},\rho_a,\mu_a)$ and $(\dot{X},\delta_a,\mu_a)$, which are $4K^2$-quasi-metric and $4K_0^2$-quasi-metric spaces (by Lemma \ref{z1}), respectively. By \eqref{q-1}, we compute  that the identity maps
$$\phi_\rho:(\dot{X},\rho)\to (\dot{X},\rho_a)\;\;\;\;\;\;\;\;\mbox{and}\;\;\;\;\;\;\;\;\phi_\delta:(\dot{X},\delta)\to (\dot{X},\delta_a)$$
are both M\"{o}bius transformations. Thus we get a quasi-m\"{o}bius correspondence
$$\varphi_a=\phi_\delta\circ \varphi \circ \phi_\rho^{-1}: (\dot{X},\rho_a)\to (\dot{X},\delta_a).$$

Now we are in a position to prove that the above identification $\varphi_a$ is quasi-symmetric with control function depending only on $C_\mu,\tau$ and $K$.

\blem \label{lem4.7}$\varphi_a: (\dot{X},\rho_a)\to (\dot{X},\delta_a)$ is quasi-symmetric. \elem
\bpf We first note from Theorem \Ref{z-4} that to prove this lemma we only need to find a constant $\lambda>0$ and a tripe of points in $X$ such that $\varphi_a$ satisfies the $\lambda$-three point condition.  Since $(X, \rho)$ is uniform perfect, by \cite[Lemma C]{WZ},  we may assume that $(X,\rho_a)$ is also $\tau$-uniform perfect.
 Since $(X,\rho)$ is a $K$-quasi-metric space, we have $1\leq \diam(\dot{X},\rho_a)\leq K$. And similarly, $1\leq \diam(\dot{X},\delta_a)\leq K_0$ because $(\dot{X},\delta_a)$ is a $K_0$-quasi-metric space. Then by Lemma \ref{z-5} we  obtain that there are constants $\alpha>0$ and $C_0>0$ depending only on $C_\mu,\tau$ and $K$ such that
$$\mu(B_\rho(a,r))\leq C_0r^\alpha$$
for all $0<r<1$ (note that we have normalized $\mu(B_\rho(a,1))=1$). Put $0<t_0<1$ satisfying
$$4K_0^2C_\mu^{\log_2(\frac{4K^3}{\tau}+1)}C_0(t_0\tau)^{\alpha}=\frac{1}{2}.$$
We point out that the number $t_0$ depends only on $C_\mu,\tau$ and $K$.

Since $(\dot{X},\rho_a)$ is $\tau$-uniformly perfect and $X\setminus B_a(a,t_0)\neq \emptyset$, there is some point $b\in B_a(a,t_0)\setminus B_a(a,\tau t_0)$  such  that  $\tau t_0\leq \rho_a(a,b)<t_0$. Since $(\dot{X},\rho_a)$ is a $4K^2$-quasi-metric space, we compute that $$1=\rho_a(a,\infty)\leq
4K^2(\rho_a(a,b)\vee \rho_a(b,\infty))=4K^2\rho_a(b,\infty).$$
This together with the choice of $t_0$ deduce
\be\label{zz5} \rho_a(a,\infty)\wedge \rho_a(b,\infty)\wedge \rho_a(a,b) \geq \tau t_0\geq \frac{\tau t_0}{K_0}\diam(\dot{X},\rho_a). \ee
On the other hand, we find
\be\label{q-3.2}\tau t_0 =\frac{\tau t_0}{1\times 1}\le \rho_a(a,b)<\rho(a,b)=\frac{\rho_a(a,b)}{\rho_a(a,\infty)\rho_a(b,\infty)}<\frac{t_0}{1\times \frac{1}{4K^2}}=4K^2t_0,\ee
which combines the estimate  $$1=\mu(B_\rho(a,1))\leq C_\mu^{\log_2(\frac{1}{\tau t_0}+1)}\mu(B_\rho(a,\tau t_0))$$ and the choice of $t_0,$  show that
\beq\label{q-3.3}
t_1 &\leq& \mu(B_\rho(a,\tau t_0))\leq \mu(B_\rho(a,\rho(a,b)\cup B_\rho(b,\rho(a,b))))
\\ \nonumber&=& \delta(a,b)\leq \mu(B_\rho(a,K\rho(a,b)))\leq C_\mu^{\log_2(\frac{4K^3}{\tau }+1)}\mu(B_\rho(a,\tau t_0))
\\ \nonumber&\leq& C_\mu^{\log_2(\frac{4K^3}{\tau }+1)}C_0(t_0\tau)^{\alpha}=\frac{1}{8K_0^2},
\eeq
where $t_1=C_\mu^{-\log_2(\frac{1}{\tau t_0}+1)}$. Moreover, by means of (\ref{q-3.2}) and (\ref{q-3.3}) we get $\delta_a(a,\infty)=1$ and $\delta_a(a,b)=\frac{\delta(a,b)}{1+\delta(a,b)}\in [\frac{t_1}{1+t_1},\frac{1}{8K_0^2}]$. This yields
$$1=\delta_a(a,\infty)<4K_0^2(\delta_a(a,b)\vee \delta_a(b,\infty))<4K_0^2(\frac{1}{8K_0^2}\vee \delta_a(b,\infty))=4K_0^2\delta_a(b,\infty).$$
Therefore
\be\label{q-3.4}\delta_a(a,\infty)\wedge \delta_a(a,b)\wedge \delta_a(a,b)\geq \frac{t_1}{4K_0^2(1+t_1)}\geq \frac{t_1}{4K_0^3(1+t_1)}\diam(\dot{X},\delta_a).\ee

Consequently, we see from (\ref{zz5}) and (\ref{q-3.4}) that $\varphi_a$ satisfies the $\lambda$-three point condition for points $a, b$ and $\infty$ with $\lambda=\frac{t_1}{4K_0^3(1+t_1)}$ depending only on $C_\mu,\tau$ and $K$. Hence this lemma follows.
\epf

Combining Theorems \ref{z-1} and \ref{z-3} with Lemma  \ref{lem4.7}, we get Theorem \ref{z-6}.
\section{The proof of main results}

In this section, we shall complete the proof of Theorems \ref{m-thm-1}, \ref{m-thm-2} and \ref{m-thm-3}.
\subsection{The proof of Theorem \ref{m-thm-1}} Since the sufficiency follows from \cite[Theorem 5.1]{WZ}, we only need to verify the necessity. To this end, we establish the following result, from which the necessity follows.

\begin{lem}\label{z-11}
A quasi-metric space is quasi-m\"{o}bius equivalent to a doubling uniformly perfect quasi-metric space, then it is also doubling and uniformly perfect.
\end{lem}

\bpf Suppose that a homeomorphism $f:(X,\rho)\to (Y,\sigma)$ is $\theta$-quasim\"{o}bius between two quasi-metric spaces, and that $Y$ is a doubling uniformly perfect space, then we need to show $X$ is also uniformly perfect and doubling as a quasi-metric space. Toward this end, according to \cite[$4.1.14$]{HKST} and \cite[Theorem $3.19$]{Vai-5}, we may assume that $X$ and $Y$ both are also complete.

We only discuss the situation whenever $(X,\rho)$ and $(Y,\sigma)$ are both bounded. For the unbounded cases, the discussion is similar. Let $c\in X$ and $c'=f(c)\in Y$. First, by means of Theorem \Ref{z-8} we know that there is a doubling measure $\nu$ defined on $(Y,\sigma)$. Then according to Lemma \ref{z1}, the flattening transformation with respect to $c\in X$ inducing an unbounded quasi-metric space $(X^c,\rho^c)$ and the identity map $\varphi^c:(X,\rho)\to (\dot{X^c},\rho^c)$ is M\"{o}bius with $\varphi^c(c)=\infty$. Also from Theorem \ref{z-6} it follows that the flattening space $(Y^{c'},\sigma^{c'})$ is an unbounded quasi-metric space which admits a doubling measure $\nu^{c'}$ and the identity map $\psi:(Y,\sigma)\to (\dot{Y^{c'}},\sigma^{c'})$ is M\"{o}bius with $\psi(c')=\infty$. Therefore, we obtain a quasi-m\"{o}bius mapping induced by $f$,
$$\widetilde{f}=\psi\circ f\circ (\varphi^{c})^{-1}: (\dot{X^c},\rho^c) \to (\dot{Y^{c'}},\sigma^{c'})$$
with $\widetilde{f}(\infty)=\infty$, which implies that $\widetilde{f}$ is quasi-symmetric by means of \cite[Theorem 3.20]{Vai-5}.

Next, invoking this fact and Lemma \ref{z-9} we see that the pull-back measure $\mu_{\widetilde{f}}$ of $\nu^{c'}$ via $\widetilde{f}$ is a doubling measure on $(\dot{X^c},\rho^c)$. Moreover, appealing to Theorem \ref{z-6} we find that the  sphericalized space $(\dot{X^c},(\rho^c)_\infty,(\mu_{\widetilde{f}})_\infty)$ is a doubling quasi-metric space. Since every doubling measure space is also doubling as a quasi-metric space, $(\dot{X^c},(\rho^c)_\infty)$ is doubling.

Furthermore, we see from Lemma \ref{z-10} that the induced map
$$\varphi=\varphi_\infty\circ \varphi^c:(X,\rho) \to (\dot{X^c},(\rho^c)_\infty)$$
is bilipschitz. Since the doubling property is clearly a bilipschitz invariant, we obtain the required assertion that $(X,\rho)$ is doubling. Hence this deduces Lemma \ref{z-11}.
\epf


\subsection{The proof of Theorem \ref{m-thm-2}} We begin the proof of Theorem \ref{m-thm-2} by showing that the sphericalization of an unbounded quasi-metric space preserves  the $Q$-Loewner property.
\begin{lem}\label{z-16}
Suppose that $(X,d,\mu)$ is an unbounded $Q$-regular $Q$-Loewner locally compact metric measure space with $Q>1$ and $a\in X$, then the  sphericalized space $(X,\widehat{d_a},\mu_a)$ is a bounded $Q$-regular $Q$-Loewner metric measure space, quantitatively.
\end{lem}
\bpf First, we see from the metric version of Theorem \ref{z-1} (see also \cite[Proposition 3.1]{LS}) that the  sphericalized space $(X,\widehat{d_a},\mu_a)$ is $Q$-regular, where $\widehat{d}_a$ is a metric on $X$ which satisfies $\frac{1}{4}d_a(x,y)\leq \widehat{d}_a(x,y)\leq d_a(x,y)$ for all $x,y\in X$ and
$$d_a(x,y)=d_a(y,x)=\begin{cases}
\displaystyle\; \frac{d(x,y)}{[1+d(x,a)][1+d(y,a)]},\;
\;\;\;\;\mbox{if}\;\;x,y\neq \infty,\\
\displaystyle\;\;\;\;\;\frac{1}{1+d(x,a)},\;\;\;\;
\;\;\;\;\;\;\;\;\;\;\;\;\;\;\;\;\; \mbox{if}\;\;y=\infty\neq x,\\
\displaystyle\;\;\;\;\;\;\;\;\;\;0,
\;\;\;\;\;\;\;\;\;\;\;\;\;\;\;\;\;\;\;\;\;\;\;\;\;\;\;\;\;\;\;\mbox{if}\;\; x=\infty=y.
\end{cases}
$$

Without loss of generality, we may assume that
$$\mu_a(A)=\int_A \rho_a^Q d\mu,$$
where $\rho_a(z)=\frac{1}{[1+d(a,z)]^2}$ and $A\subset X$ is a Borel set. Moreover, by \cite[(2.14)]{LS} we know that for any rectifiable curve $\gamma$ in $X$,
$$\ell_{\widehat{d_a}}(\gamma)=\int_\gamma \rho_a(z)|dz|.$$
From the above we deduce that $\rho_a$ is a conformal density, that is,
$$mod_Q(\Gamma,d,\mu)=mod_Q(\Gamma,d_a,\mu_a)$$
for all rectifiable family of curves $\Gamma$ in $X$. Indeed, we know that for all nonnegative Borel function $\rho:X\to[0,\infty]$ and any rectifiable curve $\gamma\in \Gamma$,
$$\int_X\rho^Qd\mu=\int_X(\frac{\rho}{\rho_a})^Qd\mu_a\;\;\;\;\mbox{and}\;\;\;\int_\gamma\rho ds=\int_\gamma \frac{\rho}{\rho_a} ds_a.$$
Therefore, according to the definition of Loewner space, we only need to find a lower bound of $\it{mod}_Q(\Gamma,d_a,\mu_a)$. It suffices to find some function $\psi:(0,\infty)\to (0,\infty)$ 
 given
$$\Delta_d(E,F)=\frac{\dist_d(E,F)}{\diam_dE\wedge \diam_dF},\;\;\;\;\;\Delta_{\widehat{d_a}}(E,F)=\frac{\dist_{\widehat{d_a}}(E,F)}{\diam_{\widehat{d_a}}E\wedge \diam_{\widehat{d_a}}F},$$
where $E$ and $F$ disjoint, nondegenerate continua in $X$, then we have the following
$$\Delta_d(E,F){\ge} \psi(\Delta_{\widehat{d_a}}(E,F)).$$

To this end, take $x\in E$ and $y\in F$ such that $\dist_d(E,F)= d(x,y)$. No loss of generality, we may assume that $\diam_dE\leq \diam_dF$. Moreover, choose $z\in E$ such that $\diam_dE\leq 2d(x,z)$ and choose  $w\in F$ such that $\diam_dF\leq 2d(y,w)$. On the other hand, since the identity map $id:(X,d)\to (X,\widehat{d_a})$ is $\theta$-quasim\"{o}bius with $\theta(t)=16t$, from \cite[Lemma 3.3]{BK} it follows that
\beq\nonumber
\frac{\widehat{d_a}(x,y)\wedge \widehat{d_a}(z,w)}{\widehat{d_a}(x,z)\wedge \widehat{d_a}(y,w)} &\leq& \theta_0(\frac{\widehat{d_a}(x,y) \widehat{d_a}(z,w)}{\widehat{d_a}(x,z) \widehat{d_a}(y,w)})
\\ \nonumber&\leq& \theta_0(16\frac{d(x,y) d(z,w)}{d(x,z) d(y,w)})
\\ \nonumber&\leq& \eta(\frac{d(x,y)\wedge d(z,w)}{d(x,z)\wedge d(y,w)}),
\eeq
where $\theta_0(t)=4(t\vee \sqrt{t})$ and $\eta(t)=\theta_0(\frac{16}{\theta_0^{-1}(1/t)})$. Consequently, from the above facts we get
\beq\nonumber
\Delta_{\widehat{d_a}}(E,F)&=& \frac{\dist_{\widehat{d_a}}(E,F)}{\diam_{\widehat{d_a}}E\wedge \diam_{\widehat{d_a}}F}\leq \frac{\widehat{d_a}(x,y)\wedge \widehat{d_a}(z,w)}{\widehat{d_a}(x,z)\wedge \widehat{d_a}(y,w)}
\\ \nonumber&\leq& \eta(\frac{d(x,y)\wedge d(z,w)}{d(x,z)\wedge d(y,w)})\leq \eta(\frac{2\dist_d(E,F)}{\diam_dE\wedge \diam_dF}).
\eeq
Hence we complete the proof of Lemma \ref{z-16} by letting $\psi(t)=\frac{1}{2}\eta^{-1}(t)$.
\epf

On the other hand, we can show that the $Q$-Loewner property is preserved under the flattening of a bounded quasimetric measure space. Since the argument for this result is completely similar  to the proof of Lemma \ref{z-16}, we do not provide the proof.

\begin{lem}\label{z-17}
Suppose that $(X,d,\mu)$ is a bounded $Q$-regular $Q$-Loewner metric measure space with $Q>1$ and $c\in X$, then the flattening space $(X^c,d^c,\mu^c)$ is an unbounded $Q$-regular $Q$-Loewner metric measure space, quantitatively.
\end{lem}

Now we are going to prove Theorem \ref{m-thm-2} by means of Lemmas \ref{z-16} and \ref{z-17}.

\emph{ The proof of Theorem \ref{m-thm-2}.} We only consider the case whenever $X$ and $Y$ are both bounded; for the other cases, it is easy to deal with and the proof is rather similar. Fix $c\in X$ and $c'=f(c)\in Y$. Then it follows from Theorem \ref{z-2} (see also \cite[Proposition 3.1]{LS}) that the flattening spaces, $(X^c,d^c,\mu^c)$ and $(Y^{c'},\sigma^{c'},\nu^{c'})$, are both $Q$-regular metric space. Moreover, by Lemma \ref{z-17} we get that $(X^c,d^c,\mu^c)$ is also $Q$-Loewner. On the other hand, a direct computation gives that the identities maps
$$\varphi_X:(X,d)\to (\dot{X^c},d^c)\;\;\;\;\mbox{and}\;\;\;\;\varphi_Y:(Y,\sigma)\to (\dot{Y^{c'}},\sigma^{c'})$$
are both $16$-quasim\"{o}bius with $\varphi_X(c)=\infty$ and $\varphi_Y(c')=\infty$. Hence we obtain a quasim\"{o}bius mapping:
$$\widehat{f}=\varphi_Y\circ f\circ \varphi_X^{-1}:(\dot{X^c},d^c)\to (\dot{Y^{c'}},\sigma^{c'})$$
with $\widehat{f}(\infty)=\infty$. Thus $\widehat{f}$ is quasi-symmetric by means of \cite[Theorem 3.20]{Vai-5}. Therefore, appealing to \cite[Corollary $1.6$]{Tys} we see that $(\dot{Y^{c'}},\sigma^{c'},\nu^{c'})$ is $Q$-Loewner. Furthermore, from Lemma \ref{z-16} it follows that the  sphericalized space $(\dot{Y^{c'}},(\sigma^{c'})_\infty,(\nu^{c'})_\infty)$ of $(\dot{Y^{c'}},\sigma^{c'},\nu^{c'})$ with respect to $\infty$ via the spherical deformation $\varphi_\infty$ is also $Q$-regular and $Q$-Loewner. Since $Q$-Loewner is a bilipschitz invariant, by Lemma \ref{z-10} we find that $(Y,\sigma,\nu)$ is also $Q$-Loewner. Hence the proof of Theorem \ref{m-thm-2} is complete.

%
%
%

\subsection{The proof of Theorem \ref{m-thm-3}}
Suppose that $X$ is $\delta$-hyperbolic for some nonnegative constant $\delta$. Thanks to \cite[Lemma 2.2.2]{BuSc}, we see that there is a constant $C=C(\delta)\geq 0$ such that
\be\label{l-1} |(\xi|\eta)_o-(\xi|\eta)_{o'}|\leq d(o,o')+C,\ee
for all $o,o'\in X$ and $\xi,\eta\in \partial_\infty X$. Then for all $0<\varepsilon<\varepsilon_0(\delta)$, by (\ref{l-0.1}) and (\ref{l-1}) we know that any two Bourdon metrics $d_{o,\varepsilon}$ and $d_{o',\varepsilon}$ are bilipschitz equivalent. Since the $Q$-regularity of $Q$-dimensional Hausdorff measure $\mathcal{H}_Q$ is  bilipschitz invariant, we may assume that $d_B$ is a Bourdon metric based at $o\in X$ with parameter $\varepsilon$ and $d_H$ is a Hamenst\"adt metric based at a Busemann function $b=b_{\omega,o}$ with parameter $\varepsilon$, and $\omega\in \partial_\infty X$. Denote
$$\rho_B(\xi,\eta)=e^{-\varepsilon (\xi|\eta)_o}\;\;\;\;\;\;\;\;\;\;\mbox{and}\;\;\;\;\;\;\;\;\;
\rho_H(\xi,\eta)=e^{-\varepsilon (\xi|\eta)_{\omega,o}},$$
with $(\xi|\eta)_{\omega,o}=(\xi|\eta)_o-(\xi|\omega)_o-(\eta|\omega)_o$. Then by \cite[Lemma 2.2.2]{BuSc} and (\ref{l-0.1}) we know that $\rho_B$ is a quasimetric on $\partial_\infty X$ and bilipschitz equivalent to $d_B$. Similarly, according to \cite[(3.4)]{BuSc}, we find that $(\partial_\infty X,\rho_H)$ and $(\partial_\infty X,d_H)$ are bilipschitz equivalent. Since the $Q$-regularity of $Q$-dimensional Hausdorff measure $\mathcal{H}_Q$ is preserved under bilipschitz mapping, it follows from the above facts that $(\partial_\infty X, d_B,\mathcal{H}_Q)$ is Ahlfors $Q$-regular if and only if the space $(\partial_\infty X,\rho_B,\mathcal{H}_Q)$ is Ahlfors $Q$-regular; $(\partial_\infty X, d_H,\mathcal{H}_Q)$ is Ahlfors $Q$-regular if and only if the space $(\partial_\infty X,\rho_H,\mathcal{H}_Q)$ is Ahlfors $Q$-regular.

Moreover, a direct computation gives that
\be\label{l-2} \rho_H(\xi,\eta)=\frac{\rho_B(\xi,\eta)}{\rho_B(\xi,\omega) \rho_B(\eta,\omega)}\ee
and so $\rho_H$, associated  with the point $\omega$, is the flatting transformation of the quasimetric $d_B$ on $\partial_\infty X$.

For the necessity, assume first that $(\partial_\infty X, d_B,\mathcal{H}_Q)$ is Ahlfors $Q$-regular and so is $(\partial_\infty X, \rho_B,\mathcal{H}_Q)$, where $\mathcal{H}_Q$ is the $Q$-dimensional Hausdorff measure. Thus by Theorem \ref{z-2} we know that the flattening measure $\mu:=(\mathcal{H}_Q)^{\omega}$ of $\mathcal{H}_Q$ with respect to $\omega$ on the space $(\partial_\infty X,\rho_H)$ is $Q$-regular as well because $\rho_H=(\rho_B)^{\omega}$ by way of (\ref{l-2}). Consequently,  following from \cite[Exercise 8.11]{Hei} it  suffices to see that the $Q$-dimensional Hausdorff measure $\mathcal{H}_Q$ defined on $(\partial_\infty X,\rho_H)$ is also $Q$-regular and we are done.

Now we are in a position to prove the sufficiency.  We assume that $(\partial_\infty X,d_H,\mathcal{H}_Q)$ is Ahlfors $Q$-regular and so is $(\partial_\infty X,\rho_H,\mathcal{H}_Q)$. Consider the  sphericalized space $(\partial_\infty X,(\rho_H)_\infty)$ associated  with the infinity $\infty$. On one hand, we see from Theorem \ref{z-1} that the corresponding  sphericalized measure $(\mathcal{H}_Q)_\infty$ of $\mathcal{H}_Q$ is Ahlfors $Q$-regular on $(\partial_\infty X,(\rho_H)_\infty)$. On the other hand, by Lemma \ref{z-10}, we know that $(\partial_\infty X, \rho_B)$ is bilipschitz equivalent to $(\partial_\infty X,(\rho_H)_\infty)$. So  the $Q$-regularity of the measure $(\mathcal{H}_Q)_\infty$ on $(\partial_\infty X, \rho_B)$ follows. Then again by \cite[Exercise 8.11]{Hei}, we obtain that $(\partial_\infty X,\rho_B,\mathcal{H}_Q)$ is Ahlfors $Q$-regular.

Hence the proof of Theorem \ref{m-thm-3} is complete.


\end{document}